%% file: iPPA_mVI140801.tex
\documentclass[final]{siamltex}
\usepackage{latexsym,url,amscd,amsmath,amssymb}
\usepackage{graphics}
\usepackage{graphicx}
\usepackage{epstopdf}
\usepackage[usenames]{color}
\usepackage{comment}
\usepackage[misc]{ifsym}
\usepackage[algo2e,linesnumbered,vlined,ruled]{algorithm2e}

\newtheorem{thm}{Theorem}
\newtheorem{assumption}{Assumption}
\newtheorem{remark}{Remark}

\usepackage{longtable}
\usepackage{color}
\def\mc{\multicolumn}

\def\Gid#1#2{\langle#1,\,#2\rangle_G}
\def\Mid#1#2{\langle#1,\,#2\rangle_M}

\def\Gnorm#1{ \|#1\|_G }
\def\Mnorm#1{ \|#1\|_M }

\def\calX{ {\cal X} }
\def\calY{ {\cal Y} }

\begin{document}
\title{A general inertial proximal point method for mixed variational inequality problem}

\author{Caihua Chen\footnotemark[1] \and Shiqian Ma\footnotemark[2]  \and  Junfeng Yang\footnotemark[3]}
\renewcommand{\thefootnote}{\fnsymbol{footnote}}

\footnotetext[1]{International Center of Management Science and Engineering, School
of Management and Engineering, Nanjing University, China (Email: {\tt chchen@nju.edu.cn}). Research supported in part
by Natural Science Foundation of Jiangsu Province under project grant No. BK20130550
and  the Natural Science Foundation of China NSFC grant 11371192.
}

\footnotetext[2]{Department of Systems Engineering and Engineering Management, The Chinese University of Hong Kong, Hong Kong (Email: {\tt sqma@se.cuhk.edu.hk}). This author was supported by Hong Kong Research Grants Council General Research Fund Early Career Scheme (Project ID: CUHK 439513).}

\footnotetext[3]{Corresponding author (Email: {\tt jfyang@nju.edu.cn}). Department of Mathematics, Nanjing University, China.
This author was supported by the Natural Science Foundation of China NSFC-11371192 and a grant form Jiangsu Key Laboratory for Numerical Simulation of Large Scale Complex Systems.
The work was done while this author was visiting the Chinese University of Hong Kong.}

\renewcommand{\thefootnote}{\arabic{footnote}}
\date{\today}
\maketitle

\begin{abstract}
In this paper, we first propose a general inertial proximal point method for the mixed variational inequality (VI) problem.
Based on our knowledge, without stronger assumptions, convergence rate result is not known in the literature for   inertial type proximal point methods.
Under certain conditions, we are able to establish the global convergence and a $o(1/k)$  convergence rate result (under certain measure) of the proposed general inertial proximal point method.
We then show that the linearized alternating direction method of multipliers (ADMM) for separable convex optimization with linear constraints is an application of a general proximal point method, provided that the algorithmic parameters are properly chosen.
As byproducts of this finding, we establish global convergence and $O(1/k)$ convergence rate results of the linearized ADMM  in both ergodic and nonergodic sense.
In particular, by applying the proposed inertial proximal point method for mixed VI to linearly constrained separable convex optimization,  we obtain an inertial version of the linearized ADMM for which the global convergence is guaranteed.
We also demonstrate the effect of the inertial extrapolation step via  experimental results on the compressive principal component pursuit problem.
\end{abstract}


\thispagestyle{plain}

\section{Introduction}
Let $T: \Re^n \rightrightarrows \Re^n$ be a set-valued maximal monotone operator from $\Re^n$ to its power set. The maximal monotone operator inclusion problem is to find $w^*\in\Re^n$ such that
\begin{equation}\label{prob:0inTw}
  0 \in T(w^*).
\end{equation}
Due to the mathematical generality of maximal monotone operators, the problem \eqref{prob:0inTw} is very inclusive and  serves as a unified model for many problems of fundamental importance, for example, fixed point problem, variational inequality problem, minimization of closed proper  convex functions, and their  extensions. Therefore, it becomes extremely important in many cases to solve \eqref{prob:0inTw} in practical and efficient ways.

The classical proximal point method, which converts the maximal monotone operator inclusion problem to a fixed point problem of a firmly nonexpansive mapping via resolvent operators, is one of the most influential approaches for solving \eqref{prob:0inTw} and has been  studied extensively both in theory and in practice.
The proximal point method was originally proposed by Martinet \cite{Mar70} based on the work of Moreau \cite{Mor65} and was popularized by Rockafellar \cite{Roc76a}.
It turns out that the proximal point method is a very powerful algorithmic tool and contains many well known algorithms as special cases. In particular, it was shown that the classical augmented Lagrangian method for constrained optimization \cite{Hes69, Pow69}, the Douglas-Rachford operator splitting method \cite{DR56} and the alternating direction method of multipliers (ADMM, \cite{GM75, GM76}) are all applications of the proximal point method, see  \cite{Roc76b, EB92}. Various inexact, relaxed and accelerated variants of the proximal point method were also very well studied in the literature, see, e.g., \cite{Roc76a, EB92, Gul92}.

The primary proximal point method for minimizing a differentiable function $f:\; \Re^n\rightarrow \Re$ can be interpreted as an implicit one-step discretization method for the ordinary differential equations
\begin{equation}\label{ode-order1-nabla-f}
  w' + \nabla f(w) = 0,
\end{equation}
where $w: \Re\rightarrow \Re^n$ is differentiable,  $w'$ denotes its derivative, and $\nabla f$ is the gradient of $f$. Suppose that $f$ is closed proper and convex and its minimum value  is attained, then every solution trajectory $\{w(t): \; t\geq 0\}$ of the differential system \eqref{ode-order1-nabla-f} converges to a minimizer of $f$ as $t$ goes to infinity. Similar conclusion can be drawn for \eqref{prob:0inTw} by considering the evolution differential inclusion problem $0 \in w'(t) + T(w(t))$ almost everywhere on $\Re_+$, provided that the operator $T$ satisfies certain conditions, see e.g., \cite{Bru75}.

The proximal point method is a one-step iterative method, i.e., each new iterate point does not depend on any iterate points already generated other than the current one.
To speed up convergence, multi-step methods have been proposed in the literature by discretizing  a second-order ordinary differential system of the form
\begin{equation}
  \label{HBF}
w'' + \gamma w' + \nabla f(w) = 0,
\end{equation}
where $\gamma>0$. Studies in this direction can be traced back to at least \cite{Pol64} which examined the system \eqref{HBF} in the context of optimization.
In the two-dimensional case, the system \eqref{HBF}  characterizes roughly the motion of a heavy ball which  rolls under its own inertial over the graph of $f$ until friction stops it at a stationary point of $f$.
The three terms in \eqref{HBF} denote, respectively, inertial force, friction force and gravity force.
Therefore, the system \eqref{HBF} is usually referred to as the \emph{heavy-ball with friction} (HBF) system.
It is easy to show that the energy function $E(t) = \frac{1}{2} \|w'(t)\|^2 + f(w(t))$ is always decreasing with time $t$ unless $w'$ vanishes, which implies that the HBF system is dissipative.
It was proved in  \cite{Alv00} that if $f$ is convex and its minimum value is attained then each solution trajectory $\{w(t): t\geq 0\}$ of \eqref{HBF} converges  to a minimizer of $f$.
In theory the convergence of the solution trajectories of the HBF system to a stationary point of $f$  can be faster than those of the first-order system \eqref{ode-order1-nabla-f}, while
in practice the second order inertial term $w''$ can be exploited to design faster algorithms \cite{APZ84,Ant94}.
Motivated by the properties of \eqref{HBF}, an implicit  discretization method was proposed in \cite{Alv00}. Specifically, given $w^{k-1}$ and $w^k$, the next point $w^{k+1}$ is determined via
\[
\frac{w^{k+1} - 2w^k + w^{k-1}}{h^2}  + \gamma\frac{w^{k+1} -  w^k}{h} +  \nabla f(w^{k+1}) = 0,
\]
which results to an iterative algorithm of the form
\begin{equation}
  \label{iPPA-min-f}
w^{k+1} = (I + \lambda \nabla f)^{-1} (w^k + \alpha (w^k-w^{k-1})),
\end{equation}
where $\lambda = h^2/(1+\gamma h)$ and $\alpha = 1/(1+\gamma h)$. Note that  \eqref{iPPA-min-f} is nothing but a proximal point step applied to the extrapolated point $w^k + \alpha (w^k-w^{k-1})$, rather than $w^k$ as in the classical proximal point method. Thus the resulting iterative scheme \eqref{iPPA-min-f}  is a two-step method and is usually referred as an inertial proximal point algorithm (PPA). Convergence properties of \eqref{iPPA-min-f} were studied in \cite{Alv00} under some assumptions on the parameters $\alpha$ and $\lambda$.
Subsequently, this inertial technique was extended to solve the inclusion problem \eqref{prob:0inTw} of maximal monotone operators in \cite{AA01}.
See also \cite{ME03} for approximate inertial PPA and \cite{Alv04, MM07, MM10} for some inertial type hybrid proximal algorithms.
Recently, there are increasing interests in studying inertial type algorithms.
Some latest references are  inertial forward-backward splitting methods for certain separable nonconvex optimization problems \cite{OCBP14}  and for  strongly convex problems \cite{OBP14, APR14}, inertial versions of the Douglas-Rachford operator splitting method and the ADMM for maximal monotone operator inclusion problem \cite{BCH14, BC14a}, and inertial  forward-backward-forward method  \cite{BC14d} based on Tseng's approach \cite{Tse00}. See also \cite{BC14b, BC14c}.

\subsection{Contributions}
In this paper, we focus on the mixed variational inequality (VI) problem and study inertial PPA under a more general setting.
In particular, a weighting matrix $G$ in the proximal term is introduced. In our setting the matrix $G$ is allowed to be positive semidefinite, as long as it is positive definite in the null space of a certain matrix.
We establish its global convergence and a $o(1/k)$ convergence rate result  under certain conditions.
To the best of our knowledge, without stronger assumptions, convergence rate result is not known in the literature for general inertial type proximal point methods.
This general setting allows us to propose an inertial version  of the linearized ADMM, a practical variant of the well-known ADMM which has recently found numerous applications \cite{Boyd+11}.
We show that the linearized ADMM for separable convex optimization is an application of a general PPA to the primal-dual optimality conditions, as long as the parameters are properly chosen. As byproducts of this finding,  we establish global convergence and $O(1/k)$ convergence rate results of the linearized ADMM.
Another aim of this paper is to study the effect of the inertial extrapolation step via numerical experiments.
Finally,  we connect inertial type algorithms with the popular accelerated methods pioneered by Nesterov \cite{Nes83} and give some concluding remarks.

The main reason that  we restrict our analysis to mixed VI problem rather than the apparently more general problem \eqref{prob:0inTw} is because it is very convenient to represent the optimality conditions of linearly constrained separable convex optimization as mixed VI. In fact, our analysis for Theorems \ref{Theorem1} and \ref{Theorem2} can be  generalized to the maximal monotone operator inclusion problem \eqref{prob:0inTw} without any difficulty.

\subsection{Notation}
We use the following notation. The standard inner product and $\ell_2$ norm are denoted by $\langle\cdot,\cdot\rangle$ and $\|\cdot\|$, respectively. The sets of symmetric, symmetric positive semidefinite and symmetric positive definite matrices of order $n$ are, respectively, denoted by $S^n, S^n_+$ and $S^n_{++}$. For any matrix $A \in S^n_{+}$ and vectors $u, v\in \Re^n$, we let $\langle u, v\rangle_A := u^TAv$ and $\|u\|_A := \sqrt{\langle u, u\rangle_A}$.  The  Frobenius norm is denoted by $\|\cdot\|_F$. The spectral radius of a square matrix $M$ is denoted by  $\rho(M)$.

\section{A general inertial PPA for mixed VI}
Let $\Omega\subseteq \Re^n$ be a closed and convex set, $\theta: \Re^n\rightarrow \Re$ be a closed proper convex function, and $F: \Re^n\rightarrow \Re^n$ be a monotone mapping.
In this paper, we consider the mixed VI problem:  find $w^*\in\Omega$ such that
\begin{equation}\label{mVI}
  \theta(w) - \theta(w^*) + \langle w-w^*,  F(w^*)\rangle \geq 0,\; \forall w\in \Omega.
\end{equation}
Let $G\in S^n_{+}$ and two sequences of parameters $\{\alpha_k\geq 0: k=0,1,2,\ldots\}$ and $\{\lambda_k > 0:   k=0,1,2,\ldots\}$ be given.
We study a general inertial PPA of the following form: given any $w^0 = w^{-1}\in \Re^n$, for $k=0,1,2,\ldots$, find $w^{k+1}\in \Omega$ such that
\begin{subequations}\label{gippa}
 \begin{eqnarray}
  \label{gippa1}
  &\hspace{-6.5cm}\bar w^k := w^k + \alpha_k (w^k - w^{k-1}),   \\
  \label{gippa2}
  &\theta(w) - \theta(w^{k+1}) +  \langle w-w^{k+1},  F(w^{k+1}) + \lambda_k^{-1}G(w^{k+1} - \bar w^k) \rangle \geq 0,   \; \forall\; w\in \Omega.
 \end{eqnarray}
\end{subequations}
We make the following assumptions.
\begin{assumption}\label{Omega-ast-nonempty}
  The set of solutions of \eqref{mVI}, denoted by $\Omega^*$, is nonempty.
\end{assumption}
\begin{assumption}\label{assumption-H-monotone}
 The mapping $F$ is  $H$-monotone in the sense that
\begin{equation}\label{H-monotonicity}
  \langle u-v, F(u)-F(v)\rangle \geq \|u-v\|_H^2, \quad \forall u, v\in \Re^n,
\end{equation}
where $H \in S^n_{+}$. Note that $H=0$ if $F$ is monotone, and $H\in S^n_{++}$ if $F$ is strongly monotone.
\end{assumption}
\begin{assumption}\label{M>0}
  The sum of $G$ and $H$, denoted by $M$, is positive definite, i.e., $M := G+H  \in S^n_{++}$.
\end{assumption}

Under Assumptions \ref{assumption-H-monotone} and \ref{M>0}, it can be shown that $w^{k+1}$ is uniquely determined in \eqref{gippa2}. Therefore, the algorithm \eqref{gippa1}-\eqref{gippa2} is well defined. Clearly, the algorithm reduces to the classical PPA if $G\in S^n_{++}$ and $\alpha_k = 0$ for all $k$. It is called inertial PPA because $\alpha_k$ can be greater than $0$. We will impose conditions on $\alpha_k$ to ensure global convergence of the inertial PPA framework \eqref{gippa}.
Our convergence results are  extensions of those in \cite{AA01}.

\begin{thm}\label{Theorem1}
Assume that Assumptions \ref{Omega-ast-nonempty}, \ref{assumption-H-monotone} and \ref{M>0} hold.
Let $\{w^k\}_{k=0}^\infty\subseteq \Re^n$ conforms to Algorithm \eqref{gippa1}-\eqref{gippa2}. The parameters $\{\alpha_k, \lambda_k\}_{k=0}^\infty$ satisfy,  for all $k$,
$0\leq \alpha_k \leq \alpha$ for some $\alpha \in [0,1)$ and $\lambda_k \geq \lambda$ for some $\lambda >0$. If
\begin{equation}\label{cond-alp-k}
  \sum_{k=1}^\infty \alpha_k \Gnorm{w^k-w^{k-1}}^2 < \infty,
\end{equation}
then the sequence $\{w^k\}_{k=0}^\infty$ converges to some point in $\Omega^*$ as $k\rightarrow\infty$.
\end{thm}
\begin{proof}
 First, we show that, for any $w^* \in \Omega^*$, $\lim_{k\rightarrow \infty} \Mnorm{w^k - w^*} $ exists. As a result,
$\{w^k\}_{k=0}^\infty$ is bounded and must have a limit point. Then, we show that any limit point of $\{w^k\}_{k=0}^\infty$ must lie in $\Omega^*$. Finally, we establish the
convergence of $\{w^k\}_{k=0}^\infty$ to a point in $\Omega^*$ as $k\rightarrow\infty$.

Let $w^*\in \Omega^*$ be  arbitrarily  chosen and $k\geq 0$. It follows from setting $w=w^* \in\Omega^*$ in \eqref{gippa2} and the $H$-monotonicity \eqref{H-monotonicity} of $F$ that
\begin{eqnarray}\label{ineq1}
\nonumber
\lambda_k^{-1}\langle w^{k+1}-w^*, w^{k+1} - \bar w^k\rangle_G
&\leq&  \theta(w^*) -  \theta(w^{k+1}) - \langle w^{k+1} - w^*,  F(w^{k+1})\rangle \\ \nonumber
&\leq&  \theta(w^*) -  \theta(w^{k+1}) - \langle w^{k+1} - w^*,  F(w^*)\rangle - \|w^{k+1} - w^*\|_H^2 \\
&\leq& - \|w^{k+1} - w^*\|_H^2.
\end{eqnarray}
%
Define $\varphi_k := \Gnorm{w^k - w^*}^2$ and recall that $\bar w^k = w^k + \alpha_k(w^k-w^{k-1})$. Plug the identities
\begin{eqnarray*}
2\Gid{w^{k+1} - w^*}{w^{k+1} - w^k}   &=& \varphi_{k+1} - \varphi_{k}   +  \Gnorm{w^{k+1}-w^k}^2,  \\
2\Gid{w^{k+1} - w^*}{w^{k} - w^{k-1}} &=& \varphi_{k}   - \varphi_{k-1} +  \Gnorm{w^{k}-w^{k-1}}^2 + 2\Gid{w^{k+1} - w^k}{w^{k}-w^{k-1}},
\end{eqnarray*}
into \eqref{ineq1} and reorganize, we obtain
\begin{eqnarray} \label{ineq2.0}
\nonumber  \psi_k & :=& \varphi_{k+1} - \varphi_{k}  -\alpha_k \left(\varphi_{k} - \varphi_{k-1}\right) \\
\nonumber
&\leq& -  \Gnorm{w^{k+1}-w^k}^2 +   2 \alpha_k \Gid{w^{k+1} - w^k}{w^{k}-w^{k-1}} + \alpha_k \Gnorm{w^{k}-w^{k-1}}^2   - 2 \lambda_k\|w^{k+1} - w^*\|_H^2 \\ \nonumber
&=& -  \Gnorm{w^{k+1}-\bar w^k}^2 +    (\alpha_k^2 + \alpha_k) \Gnorm{w^{k}-w^{k-1}}^2   - 2 \lambda_k\|w^{k+1} - w^*\|_H^2 \\ \nonumber
& \leq  & -  \Gnorm{w^{k+1}- \bar w^k}^2 +  2 \alpha_k \Gnorm{w^{k}-w^{k-1}}^2   - 2 \lambda_k\|w^{k+1} - w^*\|_H^2 \\
& \leq  & -  \Gnorm{w^{k+1}- \bar w^k}^2 +  2 \alpha_k \Gnorm{w^{k}-w^{k-1}}^2,
\end{eqnarray}
where the first inequality is due to  \eqref{ineq1} and the  second follows from $0\leq \alpha_k < 1$.
%
Define
\[\theta_k := \varphi_{k} - \varphi_{k-1} \text{~~and~~} \delta_k := 2\alpha_k \Gnorm{w^{k}-w^{k-1}}^2.\] Then, the inequality \eqref{ineq2.0} implies that
$\theta_{k+1}   \leq \alpha_k \theta_k  +  \delta_k \leq \alpha  [\theta_k]_+  +  \delta_k$,
where $[t]_+ := \max\{t,0\}$ for $t\in \Re$. Therefore, we have
\begin{eqnarray}\label{ineq+7}
[\theta_{k+1}]_+ \leq \alpha  [\theta_k]_+  +  \delta_k \leq  \alpha^{k+1}  [\theta_0]_+  + \sum_{j=0}^{k} \alpha^j \delta_{k-j}.
\end{eqnarray}
Note that by our assumption $w^0=w^{-1}$. This implies that   $\theta_0 = [\theta_0]_+ =0$
and $\delta_0=0$. Therefore, it follows from \eqref{ineq+7} that
\begin{eqnarray}\label{sum-th-k+finite}
\sum_{k=0}^\infty [\theta_{k}]_+ \leq \frac{1}{1-\alpha}   \sum_{k=0}^{\infty} \delta_{k} = \frac{1}{1-\alpha}   \sum_{k=1}^{\infty} \delta_{k}   < \infty.
\end{eqnarray}
Here the second inequality is due to the assumption \eqref{cond-alp-k}.
Let $\gamma_k := \varphi_k - \sum_{j=1}^k [\theta_j]_+$. From \eqref{sum-th-k+finite} and   $\varphi_k\geq 0$, it follows that $\gamma_k$ is bounded below.
On the other hand,
\[
\gamma_{k+1} = \varphi_{k+1} - [\theta_{k+1}]_+ - \sum_{j=1}^k [\theta_j]_+ \leq \varphi_{k+1} - \theta_{k+1} - \sum_{j=1}^k [\theta_j]_+ = \varphi_{k} - \sum_{j=1}^k [\theta_j]_+ = \gamma_k,
\]
i.e., $\gamma_k$ is nonincreasing. As a result, $\{\gamma_k\}_{k=0}^\infty$ converges as $k\rightarrow\infty$, and the following limit
\[
\lim_{k\rightarrow\infty}\varphi_k = \lim_{k\rightarrow\infty} \left(\gamma_k + \sum_{j=1}^k [\theta_j]_+\right) = \lim_{k\rightarrow\infty}  \gamma_k + \sum_{k=1}^\infty [\theta_k]_+
\]
exists. That is, $\lim_{k\rightarrow \infty} \Gnorm{w^k - w^*} $ exists for any $w^*\in \Omega^*$.
Furthermore, it follows from the second ``$\leq$" of \eqref{ineq2.0} and the definition  of $\theta_k$ and $\delta_k$ that
\begin{eqnarray}\label{ineq+4} \nonumber
  \Gnorm{w^{k+1}-  \bar w^k}^2 + 2\lambda_k\|w^{k+1}-w^*\|_H^2
&\leq& \varphi_{k} - \varphi_{k+1}  + \alpha_k \left(\varphi_{k} - \varphi_{k-1}\right)  + \delta_k   \\
&\leq& \varphi_{k} - \varphi_{k+1}  + \alpha [\theta_k]_+  + \delta_k.
\end{eqnarray}
By taking sum over $k$ and noting that $\varphi_k \geq 0$,  we obtain
\begin{equation}
  \label{ch-add1}
\sum_{k=1}^\infty \left(\|w^{k+1}-\bar{w}^k\|_G^2  + 2\lambda_k\|w^{k+1}-w^*\|_H^2 \right) \leq    \varphi_{1}  + \sum_{k=1}^\infty \left(\alpha [\theta_k]_+  + \delta_k\right) < \infty,
\end{equation}
where the second inequality follows from \eqref{sum-th-k+finite} and assumption \eqref{cond-alp-k}. Since $\lambda_k\geq \lambda >0$ for all $k$, it follows from \eqref{ch-add1} that
\begin{equation}
  \label{ch-add2}
  \lim_{k\rightarrow \infty}\|w^k-w^*\|_H = 0.
\end{equation}
Recall that $M= G+H$. Thus, $\lim_{k\rightarrow \infty} \Mnorm{w^k - w^*}$ exists.
Since $M$ is positive definite, it follows that $\{w^k\}_{k=0}^\infty$ is bounded and must have at least one limit point.

Again from \eqref{ch-add1} we have
\begin{eqnarray*}
   \lim_{k\rightarrow\infty}\Gnorm{w^{k+1}-  \bar w^k} = 0.
\end{eqnarray*}
Thus, the positive semidefiniteness of $G$ implies that $\lim_{k\rightarrow \infty} G(w^{k+1}-  \bar w^k) = 0$. On the other hand, for any fixed $w\in \Omega$,
it follows from \eqref{gippa2}  that
\begin{eqnarray}\label{gippa2-k}
  \theta(w) - \theta(w^{k}) + \langle w-w^{k},  F(w^{k})  \rangle \geq \lambda_{k-1}^{-1}\langle w^{k}-w, G(w^{k} - \bar w^{k-1})\rangle.
\end{eqnarray}
Suppose that $w^\star$ is any limit point of $\{w^k\}_{k=0}^\infty$ and $w^{k_j} \rightarrow w^\star$ as $j\rightarrow \infty$. Since $\Omega$ is closed, $w^\star\in \Omega$. Furthermore, by taking the limit over $k=k_j \rightarrow \infty$ in \eqref{gippa2-k} and noting that $G(w^{k} - \bar w^{k-1})\rightarrow 0$ and $\lambda_{k-1} \geq \lambda >0$, we obtain
\[
\theta(w) - \theta(w^\star) + \langle w-w^\star,  F(w^\star)  \rangle \geq 0.
\]
Since $w$ can vary arbitrarily in $\Omega$, we conclude that $w^\star\in \Omega^*$. That is, any limit point of $\{w^k\}_{k=0}^\infty$ must also lie in $\Omega^*$.

Finally, we establish the uniqueness of limit points of $\{w^k\}_{k=0}^\infty$. Suppose that $w^*_1$ and $w^*_2$ are two limit points of $\{w^k\}_{k=0}^\infty$ and
$\lim_{j\rightarrow \infty} w^{i_j} = w^*_1$, $\lim_{j\rightarrow \infty} w^{k_j} = w^*_2$. Assume that $\lim_{k\rightarrow \infty} \|w^k - w^*_i\|_M = v_i$ for $i=1,2$.
By taking the limit over $k=i_j\rightarrow \infty$ and $k=k_j\rightarrow \infty$ in the equality
\[
\Mnorm{w^k - w^*_1}^2 - \Mnorm{w^k - w^*_2}^2 = \Mnorm{w^*_1 - w^*_2}^2 + 2\Mid{w^*_1 - w^*_2}{w^*_2 - w^k},
\]
we obtain  $v_1 - v_2 = - \Mnorm{w^*_1 - w^*_2}^2 = \Mnorm{w^*_1 - w^*_2}^2$. Thus, $\Mnorm{w^*_1 - w^*_2} = 0$.  Since $M$ is positive definite, this implies that $w^*_1 = w^*_2$. Therefore, $\{w^k\}_{k=0}^\infty$ converges to some point in $\Omega^*$ and the proof of the theorem is completed.
\end{proof}

We have the following remarks on the assumptions and results of Theorem \ref{Theorem1}.
\begin{remark}
  In practice, it is not hard to  select $\alpha_k$ online such that the condition \eqref{cond-alp-k} is satisfied.
\end{remark}
\begin{remark}
If $\alpha_k=0$ for all $k$, then the condition \eqref{cond-alp-k} is obviously satisfied. In this case, we reestablished the convergence of the classical PPA under the weaker condition
that $G \in S^n_{+}$, provided that $\lambda_k\geq \lambda >0$ and $H+G\in S^n_{++}$, e.g., when $F$ is strongly monotone, i.e., $H\in S^n_{++}$.
\end{remark}
\begin{remark}
   Suppose that $H=0$ and $G\in S^n_{+}$, but $G \notin S^n_{++}$. Then, the sequence $\{w^k\}_{k=0}^\infty$ may not be well defined since \eqref{gippa2} does not necessarily have a solution in general. In the case that $\{w^k\}_{k=0}^\infty$ is indeed well defined (which is possible), the conclusion that $\lim_{k\rightarrow \infty} \Gnorm{w^k - w^*} $ exists for any $w^*\in \Omega^*$ still holds under condition \eqref{cond-alp-k}. However, since $G$ is  only positive semidefinite,   the boundedness of $\{w^k\}_{k=0}^\infty$ cannot be guaranteed.
If a limit point $w^\star$ of $\{w^k\}_{k=0}^\infty$ does exist, then   the conclusion $w^\star \in \Omega^*$  holds still. Moreover, suppose that $w^\star_1$ and $w^\star_2$ are any two limit points of $\{w^k\}_{k=0}^\infty$, then it holds that  $Gw^\star_1=Gw^\star_2$.
\end{remark}

In the following theorem, we remove the assumption \eqref{cond-alp-k} by assuming that the sequence $\{\alpha_k\}_{k=0}^\infty$ satisfies some additional easily implementable conditions.
Moreover, we establish a $o(1/k)$ convergence rate result for the general inertial proximal point method \eqref{gippa}. The trick used here to improve convergence rate from $O(1/k)$ to $o(1/k)$
seems to be first introduced in \cite{DLY13,DY14}.
To the best of our knowledge,  there is no  convergence rate result known in the literature without stronger assumptions for inertial type proximal point methods.
\begin{thm}\label{Theorem2}
Assume that Assumptions \ref{Omega-ast-nonempty}, \ref{assumption-H-monotone} and \ref{M>0} hold.
Suppose that the parameters $\{\alpha_k, \lambda_k\}_{k=0}^\infty$ satisfy, for all $k$, $0\leq \alpha_k \leq \alpha_{k+1}\leq \alpha < \frac{1}{3}$ and $\lambda_k\geq \lambda$ for some $\lambda>0$. Let $\{w^k\}_{k=0}^\infty$ be   the sequence  generated by  Algorithm \eqref{gippa1}-\eqref{gippa2}. Then, we have the following results.
\begin{enumerate}
  \item  $\{w^k\}_{k=0}^\infty$ converges to some point in $\Omega^*$ as $k\rightarrow\infty$;
  \item For any $w^*\in\Omega^*$ and positive integer $k$, it holds that
\begin{equation}
  \label{ineq+2}
  \min_{0\leq i\leq k-1}  \|w^{i+1}-  \bar w^i\|_G^2 \leq
\frac{\left(1 + \frac{2}{1-3\alpha}\right)\Gnorm{w^0 - w^*}^2}{k}.
\end{equation}
Moreover, it holds as $k\rightarrow\infty$ that
\begin{equation}
  \label{ineq+6}
  \min_{0\leq i\leq k-1}  \|w^{i+1}-  \bar w^i\|_G^2 = o\left(\frac{1}{k}\right).
\end{equation}
\end{enumerate}
\end{thm}
\begin{proof}
Let $w^*\in \Omega^*$ be arbitrary fixed and,  for all $k\geq 0$, retain the notation $\varphi_k = \Gnorm{w^k - w^*}^2$,
\[
\psi_k = \varphi_{k+1} - \varphi_{k}  -\alpha_k \left(\varphi_{k} - \varphi_{k-1}\right) \text{~~and~~} \theta_k = \varphi_k - \varphi_{k-1}.
\]
It follows from the first ``$\leq$'' in  \eqref{ineq2.0}  and  $\lambda_k\geq 0$ that 
\begin{eqnarray} \label{ineq5} \nonumber
\psi_k &\leq&  -  \Gnorm{w^{k+1}-  w^k}^2 +    2\alpha_k \Gid{w^{k+1}-  w^k}{w^{k}-  w^{k-1}}  +   \alpha_k \Gnorm{w^{k}-w^{k-1}}^2 \\ \nonumber
       &\leq&  -  \Gnorm{w^{k+1}-  w^k}^2 +     \alpha_k \left(\Gnorm{w^{k+1}-  w^k}^2 + \Gnorm{w^{k}-  w^{k-1}}^2\right)  +   \alpha_k \Gnorm{w^{k}-w^{k-1}}^2 \\
       &=& - (1-\alpha_k) \|w^{k+1}-  w^k\|_G^2 +  2\alpha_k \Gnorm{w^{k}-  w^{k-1}}^2,
\end{eqnarray}
where the second ``$\leq$" follows from the Cauchy-Schwartz inequality.
Define
\[ \mu_k := \varphi_k - \alpha_k \varphi_{k-1} + 2\alpha_k \Gnorm{w^{k}-  w^{k-1}}^2.\]
From $0\leq \alpha_k \leq \alpha_{k+1} \leq \alpha < \frac{1}{3}$, the fact that $\varphi_k\geq 0$ and \eqref{ineq5}, we have
\begin{eqnarray} \label{ineq7}
\nonumber    \mu_{k+1} - \mu_k
&=& \varphi_{k+1} - \alpha_{k+1} \varphi_{k} +  2\alpha_{k+1} \Gnorm{w^{k+1}-  w^{k}}^2   - \left(\varphi_k - \alpha_k \varphi_{k-1} + 2\alpha_k \Gnorm{w^{k}-  w^{k-1}}^2\right)\\ \nonumber
&\leq& \psi_k + 2\alpha_{k+1} \Gnorm{w^{k+1}-  w^{k}}^2   -  2\alpha_k \Gnorm{w^{k}-  w^{k-1}}^2\\ \nonumber
&\leq& -   (1-\alpha_k) \|w^{k+1}-  w^k\|_{G}^2 +    2 \alpha_{k+1} \Gnorm{w^{k+1}-  w^{k}}^2  \\ \nonumber
&\leq& -   (1-3\alpha  ) \|w^{k+1}-  w^k\|_{G}^2 \\
& \leq& 0.
\end{eqnarray}
Thus, $\mu_{k+1} \leq \mu_k$ for all $k\geq 0$. Note that $w^0=w^{-1}$ by our assumption. It follows from the definitions of $\mu_k$ and $\varphi_k$ that $\mu_0 = (1-\alpha_0) \varphi_0 \leq \varphi_0 := \Gnorm{w^0-w^*}^2$. Therefore, we have
\begin{eqnarray}\label{ineq+1}
- \alpha \varphi_{k-1}  \leq \varphi_k - \alpha  \varphi_{k-1} \leq \varphi_k - \alpha_k  \varphi_{k-1}  \leq \mu_k \leq \mu_0 \leq \varphi_0.
\end{eqnarray}
Further take into account \eqref{ineq7}, we obtain
\begin{equation}
   \label{ineq8} \varphi_k \leq  \alpha  \varphi_{k-1} + \varphi_0 \leq \alpha^k \varphi_0 + \varphi_0 \sum_{j=0}^{k-1}\alpha^j \leq \alpha^k \varphi_0 + \frac{\varphi_0}{1-\alpha}.
\end{equation}
The second last ``$\leq$" in \eqref{ineq7} implies that $(1-3\alpha  ) \|w^{k+1}-  w^k\|_G^2 \leq \mu_k - \mu_{k+1}$ for $k\geq 0$.
Together with  \eqref{ineq+1} and \eqref{ineq8}, this  implies
\begin{equation}
  \label{ineq+3}
(1-3\alpha  )  \sum_{j=0}^k \|w^{j+1}-  w^j\|_G^2 \leq   \mu_0 - \mu_{k+1} \leq \varphi_0 + \alpha \varphi_k \leq \alpha^{k+1}\varphi_0 + \frac{\varphi_0}{1-\alpha} \leq  2\varphi_0,
\end{equation}
where the second inequality is due to $\mu_0\leq \varphi_0$ and $- \alpha \varphi_{k} \leq \mu_{k+1}$,  the next one follows from \eqref{ineq8}, and the last one is  due to  $\alpha < 1/3$.
By taking the limit $k\rightarrow \infty$, we obtain
\begin{equation}\label{def:C1}
 \frac{1}{2}\sum_{k=1}^\infty \delta_k  = \sum_{k=1}^\infty \alpha_k\|w^{k}-  w^{k-1}\|_G^2 \leq \alpha\sum_{k=1}^\infty \|w^{k}-  w^{k-1}\|_G^2   \leq   \frac{2\varphi_0\alpha}{ 1-3\alpha } := C_1 < \infty.
\end{equation}
The convergence of $\{w^k\}_{k=0}^\infty$ to a solution point in $\Omega^*$ follows from the proof of Theorem \ref{Theorem1}.

It follows from \eqref{ineq+4} that, for $i\geq 0$,
$  \Gnorm{w^{i+1}-  \bar w^i}^2 \leq \varphi_{i} - \varphi_{i+1}  + \alpha [\theta_i]_+  + \delta_i $, from which we obtain
\begin{eqnarray}\label{ineq+5}
\sum_{i=0}^{k-1} \Gnorm{w^{i+1}-  \bar w^i}^2
 \leq \varphi_{0} - \varphi_{k}  + \alpha \sum_{i=1}^{k-1} [\theta_i]_+  + \sum_{i=1}^{k-1}  \delta_i \leq
\varphi_{0}  + \alpha C_2  + 2C_1,
\end{eqnarray}
where $C_1$ is defined in \eqref{def:C1} and $C_2$ is defined as
\[
 C_2 := \frac{   2C_1}{1-\alpha}   \geq  \frac{1}{1-\alpha}   \sum_{i=1}^{\infty} \delta_{i} \geq  \sum_{i=1}^{\infty} [\theta_i]_+.
\]
Here the first ``$\geq$" follows from the definition of $C_1$ in \eqref{def:C1} and the second one follows from \eqref{sum-th-k+finite}.
Direct calculation shows that
\begin{eqnarray}\label{calculateC}
  \varphi_{0}  + \alpha C_2  + 2C_1 = \left[1 + \left(\frac{   2\alpha  }{1-\alpha} + 2 \right)  \frac{2 \alpha}{ 1-3\alpha }\right]\varphi_0 \leq \left(1 + \frac{2}{1-3\alpha}\right)\varphi_0,
\end{eqnarray}
where the ``$\leq$" follows from $\alpha<1/3$.
The estimate \eqref{ineq+2} follows immediately from \eqref{ineq+5} and \eqref{calculateC}.
The $o\left(1/k\right)$ result \eqref{ineq+6} follows from
\begin{equation}\label{y-add1}
\frac{k-1}{2} \min_{0\leq i\leq k-1} \Gnorm{w^{i+1}-  \bar w^i}^2 \leq \sum_{i= \lfloor {k-1 \over 2} \rfloor}^{k-1} \Gnorm{w^{i+1}-  \bar w^i}^2,
\end{equation}
where $\lfloor {(k-1)/2} \rfloor$ denotes the greatest integer no greater than $(k-1)/2$, and the fact that
the right-hand-side of \eqref{y-add1}
converges to $0$ as $k\rightarrow\infty$ because $\sum_{i=0}^{\infty}\Gnorm{w^{i+1}-  \bar w^i}^2 < \infty$.
\end{proof}

\begin{remark}
  Note that $w^{k+1}$ is obtained via a proximal point step from $\bar w^k$. Thus, the equality $w^{k+1} = \bar w^k$ implies that $w^{k+1}$ is already a solution of \eqref{mVI} (even if $G$ is only positive semidefinite, see \eqref{gippa2}). In this sense, the error estimate given in \eqref{ineq+2} can be viewed as a convergence rate result of the general inertial proximal point method \eqref{gippa}. In particular, \eqref{ineq+2} implies that, to obtain an $\varepsilon$-optimal solution in the sense that $\Gnorm{w^{k+1} -\bar w^k}^2 \leq \varepsilon$, the upper bound of iterations required by \eqref{gippa} is
\[
\frac{\left(1 + \frac{2}{1-3\alpha}\right)\Gnorm{w^0 - w^*}^2}{\varepsilon}.
\]
\end{remark}
\begin{remark}
  In general Hilbert space, weak convergence of $\{w^k\}_{k=0}^\infty$ to a point in $\Omega^*$  can still be guaranteed under similar assumptions. The analysis is similar to that of Theorems \ref{Theorem1} and \ref{Theorem2} by using a well-known result, called Opial's lemma \cite{Opi67}, in functional analysis of Banach space.
\end{remark}

\section{Inertial linearized ADMM}

In this section, we prove that under suitable conditions the linearized ADMM is an application of PPA with weighting matrix $G \in S^n_{++}$.
As byproducts of this result, we establish convergence, ergodic and nonergodic convergence rate results for linearized ADMM within the PPA framework.
Furthermore, an inertial version of the linearized ADMM is proposed, whose convergence is guaranteed by Theorems \ref{Theorem1} and \ref{Theorem2}.

Let $f: \Re^{n_1}\rightarrow\Re$ and $g: \Re^{n_2}\rightarrow\Re$ be closed  convex functions, ${\cal X}\subseteq \Re^{n_1}$ and ${\cal Y}\subseteq \Re^{n_2}$ be closed convex sets. Consider linearly constrained separable convex optimization problem of the form
\begin{equation}\label{min-(f+g)}
\min_{x,y} \left\{ f(x) + g(y): \ \hbox{s.t. } Ax + By = b, x\in {\cal X}, y\in {\cal Y}\right\},
\end{equation}
where $A\in \Re^{m\times n_1}$, $B\in \Re^{m\times n_2}$ and $b\in \Re^m$ are given. We assume that the set of  KKT points  of \eqref{min-(f+g)} is nonempty.
Under very little assumptions, see, e.g., \cite{Eck89}, \eqref{min-(f+g)} is equivalent to the mixed variational inequality problem \eqref{mVI}
with  $\Omega$, $w$, $\theta$ and $F$ given, respectively, by $\Omega := {\cal X} \times {\cal Y} \times \Re^m$,
\begin{equation}\label{def-w-F-theta}
w = \left(
      \begin{array}{c}
        x \\
        y \\
        p \\
      \end{array}
    \right), \quad
\theta(w) := f(x) + g(y), \quad
F(w) =
\left(
  \begin{array}{ccc}
    0 & 0  & -A^T \\
    0 & 0  & -B^T \\
    A & B  & 0 \\
  \end{array}
\right)
\left(
      \begin{array}{c}
        x \\
        y \\
        p \\
      \end{array}
    \right)
-
\left(
  \begin{array}{c}
    0 \\
    0 \\
    b \\
  \end{array}
\right).
\end{equation}
Since the coefficient matrix defining $F$ is skew-symmetric, $F$ is monotone, and thus Assumption  \ref{assumption-H-monotone} is satisfied with $H=0$.
Let $\beta >0$ and define the Lagrangian and the augmented Lagrangian functions, respectively, as
\begin{subequations}\label{def:L-AL}
\begin{eqnarray}\label{def:L}
  {\cal L}(x,y,p) &:=& f(x) + g(y) - \langle p, Ax + By - b\rangle,\\
  \label{def:AL}
  {\cal \bar L}(x,y,p) &:=& {\cal L}(x,y,p) + \frac{\beta}{2}\|Ax+By-b\|^2.
\end{eqnarray}
\end{subequations}
Given $(y^k,p^k)$, the classical ADMM in ``$x-p-y$" order iterates as follows:
\begin{subequations}\label{ADM-xpy}
\begin{eqnarray}
\label{ADM-xpy-x}
  x^{k+1} &=& \arg\min_{x\in {\cal X}} {\cal \bar L}(x,y^k,p^k), \\
\label{ADM-xpy-y}
  p^{k+1} &=& p^k - \beta (Ax^{k+1}+By^k-b), \\
\label{ADM-xpy-y}
  y^{k+1} &=& \arg\min_{y\in {\cal Y}} {\cal \bar L}(x^{k+1},y,p^{k+1}).
\end{eqnarray}
\end{subequations}
Note that here we still use  the latest value of each variable in each step of the alternating computation. Therefore, it is equivalent to the commonly seen ADMM
in ``$y-x-p$" order in a cyclic sense. We use the order ``$x-p-y$" because the resulting algorithm can be easily explained as a PPA-like algorithm applied to the primal-dual optimality conditions, see \cite{CGHY13}.

Given $(x^k,y^k,p^k)$ and two parameters $\tau, \eta>0$, the iteration of linearized ADMM in ``$x-p-y$" order appears as
\begin{subequations}
  \label{LADM-xpy}
\begin{eqnarray}
\label{LADM-xpy-u}
  u^k     &=& A^T(Ax^k+By^k-b),\\
\label{LADM-xpy-x}
  x^{k+1} &=& \arg\min_{x\in\calX} f(x) - \langle p^k, Ax \rangle + \frac{\beta}{2\tau}\|x - (x^k -\tau u^k) \|^2, \\
\label{LADM-xpy-p}
  p^{k+1} &=& p^k - \beta (Ax^{k+1}+By^k-b), \\
\label{LADM-xpy-v}
  v^k &=& B^T(Ax^{k+1}+By^k-b),\\
\label{LADM-xpy-y}
  y^{k+1} &=& \arg\min_{y\in\calY} g(y) - \langle p^{k+1}, By\rangle + \frac{\beta}{2\eta}\|y - (y^k -\eta v^k) \|^2.
\end{eqnarray}
\end{subequations}
In the following, we prove that  under suitable assumptions $(x^{k+1},y^{k+1},p^{k+1})$ generated by \eqref{LADM-xpy} conforms to the classical PPA with an appropriate symmetric and positive definite weighting matrix $G$.
\begin{thm}\label{Theorem-LADM-mVI-k+1}
  Given $w^k = (x^k,y^k,p^k) \in \Omega$, then $w^{k+1} = (x^{k+1},y^{k+1},p^{k+1})$ generated by the linearized ADMM framework \eqref{LADM-xpy} satisfies
  \begin{equation}\label{LADM-mVI-k+1}
  w^{k+1} \in \Omega, \;     \theta(w) - \theta(w^{k+1}) + \langle w - w^{k+1}, F(w^{k+1}) + G (w^{k+1} -  w^k)\rangle \geq 0, \; \forall w\in \Omega,
  \end{equation}
 where
 \begin{equation}\label{G-LADM}
 G = \left(
       \begin{array}{ccc}
         \beta \left(\frac{1}{\tau}I - A^TA\right) & {\bf 0} & {\bf 0} \\
         {\bf 0} & \frac{\beta}{\eta} I & -B^T \\
         {\bf 0} & -B & \frac{1}{\beta}I \\
       \end{array}
     \right).
 \end{equation}
Here $I$ denotes identity matrix of appropriate size.
\end{thm}
\begin{proof}
The optimality conditions of \eqref{LADM-xpy-x} and \eqref{LADM-xpy-y} imply that
\begin{eqnarray*}
f(x) - f(x^{k+1}) + (x-x^{k+1})^T  \left\{-A^T p^k + \frac{\beta}{\tau} (x^{k+1} - x^k) + \beta A^T(Ax^k+By^k-b)\right\} \geq 0, \; \forall x\in\calX,\\
g(y) - g(y^{k+1}) + (y-y^{k+1})^T  \left\{-B^T p^{k+1} + \frac{\beta}{\eta}(y^{k+1} - y^k) + \beta B^T(Ax^{k+1}+By^k-b)\right\} \geq 0, \; \forall y\in\calY.
\end{eqnarray*}
By noting \eqref{LADM-xpy-p}, the above relations can be rewritten as
\begin{subequations}\label{LADM-op-xy}
\begin{eqnarray}
\label{LADM-op-x}
f(x) - f(x^{k+1}) + (x-x^{k+1})^T  \left\{-A^Tp^{k+1} + \beta \left(\frac{1}{\tau}I - A^TA\right) (x^{k+1} - x^k) \right\} \geq 0, \; \forall x\in\calX,\\
\label{LADM-op-y}
g(y) - g(y^{k+1}) + (y-y^{k+1})^T  \left\{   -B^T p^{k+1} + \frac{\beta}{\eta}(y^{k+1} - y^k) -   B^T (p^{k+1} - p^k)\right\} \geq 0, \; \forall y\in\calY.
\end{eqnarray}
\end{subequations}
Note that \eqref{LADM-xpy-p} can be equivalently represented as
\begin{equation}\label{LADM-op-p}
(p- p^{k+1})^T \left\{(Ax^{k+1}+By^{k+1}-b)  - B (y^{k+1} - y^k) + \frac{1}{\beta} (p^{k+1} - p^k)  \right\} \geq 0, \;\forall p\in\Re^m.
\end{equation}
By the notation defined in \eqref{def-w-F-theta}, we see that the addition of \eqref{LADM-op-x}, \eqref{LADM-op-y} and \eqref{LADM-op-p} yields \eqref{LADM-mVI-k+1}, with $G$ defined in \eqref{G-LADM}.
\end{proof}

\begin{remark}
Clearly, the matrix $G$ defined in \eqref{G-LADM} is symmetric and positive definite provided that the parameters $\tau$ and $\eta$ are reasonably small.
In particular, $G$ is positive definite if $\tau < 1/\rho(A^TA)$ and $\eta<1/\rho(B^TB)$.
Using similar analysis, it is easy to verify that $w^{k+1} = (x^{k+1},y^{k+1},p^{k+1})$ generated by the ADMM framework \eqref{ADM-xpy} conforms to
 \eqref{LADM-mVI-k+1} with $G$ defined by
  \begin{equation}\label{G-ADM}
 G = \left(
       \begin{array}{ccc}
         {\bf 0} & {\bf 0} & {\bf 0} \\
         {\bf 0} & \beta B^TB & -B^T \\
         {\bf 0} & -B & \frac{1}{\beta}I \\
       \end{array}
     \right),
 \end{equation}
which is clearly never positive  definite.  See \cite{CGHY13} for details.
 %
\end{remark}

For the linearized ADMM framework \eqref{LADM-xpy}, we have the following convergence results. Their proofs are given in the Appendix for convenience of readers.
Similar convergence analysis and complexity results can be found in \cite{HY12a,HY12b}, and also \cite{ST14}, where a unified analysis of the proximal method of multipliers is given.

\begin{thm}\label{Theorem-convergence-LADM}
Assume that $0 < \tau < 1/\rho(A^TA)$ and $0<\eta < 1/\rho(B^TB)$.
Let $\{w^k = (x^k,y^k,p^k)\}_{k=0}^\infty$ be generated by the linearized ADMM framework \eqref{LADM-xpy} from any starting point    $w^0 = (x^0,y^0,p^0)$.
The following results hold.
\begin{enumerate}
  \item The sequence $\{w^k = (x^k,y^k,p^k)\}_{k=0}^\infty$ converges to a solution of \eqref{mVI}, i.e., there exists $w^\star  = (x^\star,y^\star,p^\star) \in\Omega^*$ such that $\lim_{k\rightarrow\infty}w^k = w^\star$. Moreover, $(x^\star,y^\star) $ is a solution of \eqref{min-(f+g)}.
  \item   For any fixed integer $k>0$, define $\bar w^k := \frac{1}{k+1} \sum_{i=0}^k w^{i+1}$. Then, it holds that
\begin{equation}\label{ergodic-rate}
\bar w^k\in\Omega, \;
  \theta(w) - \theta(\bar w^k)  + ( w - \bar w^k)^T F(w)    \geq - \frac{\|w - w^0\|_G^2}{2(k+1)}, \; \forall  w\in \Omega,
\end{equation}
or, equivalently,
  \begin{equation}\label{saddle-conditions-approximate}
\bar w^k = (\bar x^k, \bar y^k, \bar p^k)\in\Omega, \; {\cal L}(\bar x^k,\bar y^k,p) -  {\cal L}(x,y,\bar p^k) \leq \frac{\|w - w^0\|_G^2}{2(k+1)}, \; \forall w = (x,y,p) \in \Omega.
  \end{equation}
Here ${\cal L}$ is the  Lagrangian  function defined in \eqref{def:L}.
  \item  After $k>0$ iterations, we have
  \begin{equation}\label{non-ergodic-rate}
    \| w^{k} - w^{k-1}\|_G^2 \leq   \frac{\|w^0 - w^*\|_G^2} {k}.
  \end{equation}
Moreover, it holds  as $k\rightarrow\infty$ that
\begin{eqnarray}\label{non-ergodic-rate-o}
  \|w^k - w^{k-1}\|_G^2 = o\left(1/k\right).
\end{eqnarray}
\end{enumerate}

\end{thm}

\begin{remark}
It is not hard to show that the set of solutions $\Omega^*$ of the mixed VI problem \eqref{mVI} can be expressed as the intersection of
$\Omega_w := \left\{  \bar w\in\Omega \ |\  \theta(w) - \theta(\bar w) + (w-\bar w)^TF(w) \geq 0   \right\}$ for all $w\in \Omega$, i.e.,
   \[\Omega^* = \bigcap_{w\in \Omega} \Omega_w =  \bigcap_{w\in \Omega} \left\{  \bar w\in\Omega \ |\  \theta(w) - \theta(\bar w) + (w-\bar w)^TF(w) \geq 0   \right\}.\]
See, e.g., \cite{FP03book}.
Therefore, the result  \eqref{ergodic-rate} essentially assures that after $k$ iterations an approximate solution $\bar w^k$ with accuracy $O(1/k)$ can be found.
On the other hand, it is easy to show that  $w^* = (x^*,y^*,p^*) \in \Omega^*$ if and only if
\begin{equation}\label{saddle-conditions}
{\cal L}(x^*,y^*,p) - {\cal L}(x,y,p^*) \leq 0, \; \forall  w = (x,y,p) \in \Omega.
\end{equation}
Thus,  \eqref{saddle-conditions-approximate} can  be viewed as an approximation to the optimality condition \eqref{saddle-conditions}.
%
%
Since $\bar w^k$ is the average of all the points generated in the first $(k+1)$ iterations, the result \eqref{ergodic-rate} or \eqref{saddle-conditions-approximate} is usually called an ergodic convergence rate.
\end{remark}
\begin{remark}
  It is easy to see from \eqref{LADM-mVI-k+1} that $w^{k+1}$ must be a solution if $w^{k+1} = w^k$. As such, the difference of two consecutive iterations can be viewed in some sense as a   measure of how close the current point is to the solution set. Therefore, the result \eqref{non-ergodic-rate} estimates the convergence rate of $w^k$ to the solution set using the measure $\|w^k-w^{k-1}\|_G^2$.
\end{remark}

\begin{remark}
  We note that all the results given in Theorem \ref{Theorem-convergence-LADM} remain valid if
the conditions  on $\tau$ and $\eta$ 
are relaxed to  $0 < \tau \leq 1/\rho(A^TA)$ and $0<\eta \leq 1/\rho(B^TB)$, respectively. The proof is a little bit complicated and we refer interested readers to \cite{FPST13,Chen12Thesis}.
\end{remark}

Now we state the inertial version of the linearized ADMM, which is new to the best of our knowledge.  Given $\beta, \tau, \eta>0$, a sequence $\{\alpha_k\geq 0\}_{k=0}^{\infty}$,
$(x^k,y^k,p^k)$ and $(x^{k-1},y^{k-1},p^{k-1})$, the inertial linearized ADMM iterates as follows:
\begin{subequations}
  \label{iLADM-xpy}
\begin{eqnarray}
\label{iLADM-xpy-bar-p}
(\bar x^k,\bar y^k,\bar p^k) &=& (x^k,y^k,p^k) + \alpha_k  (x^k-x^{k-1},y^k-y^{k-1},p^k-y^{k-1})  \\
\label{iLADM-xpy-u}
  u^k     &=& A^T(A\bar x^k+B\bar y^k-b),\\
\label{iLADM-xpy-x}
  x^{k+1} &=& \arg\min_{x\in\calX} f(x) - \langle \bar p^k, Ax \rangle + \frac{\beta}{2\tau}\|x - (\bar x^k -\tau u^k) \|^2, \\
\label{iLADM-xpy-p}
  p^{k+1} &=& \bar p^k - \beta (Ax^{k+1}+B\bar y^k-b), \\
\label{iLADM-xpy-v}
  v^k &=& B^T(Ax^{k+1}+B\bar y^k-b),\\
\label{iLADM-xpy-y}
  y^{k+1} &=& \arg\min_{y\in\calY} g(y) - \langle p^{k+1}, By\rangle + \frac{\beta}{2\eta}\|y - (\bar y^k -\eta v^k) \|^2.
\end{eqnarray}
\end{subequations}
The following convergence result is a consequence of Theorems \ref{Theorem2} and \ref{Theorem-LADM-mVI-k+1}.
\begin{thm}\label{Theorem-iLADM-xpy}
Let $G$ be defined in \eqref{G-LADM}
and $\{(x^k,y^k,p^k)\}_{k=0}^\infty\subseteq \Re^n$ be generated by \eqref{iLADM-xpy} from any starting point $(x^0,y^0,p^0)=(x^{-1},y^{-1},p^{-1})$.
Suppose  that $0<\tau<1/\rho(A^TA)$, $0<\eta<1/\rho(B^TB)$ and $\{\alpha_k\}_{k=0}^\infty$ satisfies,  for all $k$, $0\leq \alpha_k \leq \alpha_{k+1}\leq \alpha < \frac{1}{3}$. Then,  the sequence $\{(x^k,y^k,p^k)\}_{k=0}^\infty$converges to some point in $\Omega^*$, the set of solutions of \eqref{mVI}, as $k\rightarrow\infty$.
%
Moreover, it holds that
\begin{equation}
  \label{ineq+3forLADMM}
  \min_{0\leq i\leq k-1}  \|(x^{i+1},y^{i+1},p^{i+1})  - (\bar x^i,\bar y^i,\bar p^i)\|_G^2 = o\left(\frac{1}{k}\right).
\end{equation}
\end{thm}

\section{Numerical Results}
In this section, we present numerical results to compare the performance of the linearized ADMM \eqref{LADM-xpy} (abbreviated as LADMM) and the proposed inertial linearized ADMM \eqref{iLADM-xpy} (abbreviated as iLADMM).
Both algorithms were implemented in MATLAB. All the experiments were performed with Microsoft Windows 8 and MATLAB v7.13 (R2011b),
running on a 64-bit Lenovo laptop with an Intel Core i7-3667U CPU at 2.00 GHz and 8 GB of memory.

\subsection{Compressive principal component pursuit}
In our experiments, we focused on the compressive principal component pursuit problem \cite{WGMM13}, which aims to recover low-rank and sparse components from compressive or incomplete measurements. Let ${\cal A}: \Re^{m\times n} \rightarrow \Re^q$ be a linear operator, $L_0$ and $S_0$ be, respectively, low-rank and sparse matrices of size $m\times n$. The incomplete measurements are given by $b = {\cal A}(L_0+S_0)$.
Under certain technical conditions, such as $L_0$ is $\mu$-incoherent, the support of $S_0$ is randomly distributed with nonzero probability $\rho$ and the signs of $S_0$  conform  to Bernoulli distribution, it was proved in \cite{WGMM13} that the low-rank and the sparse components $L_0$ and $S_0$ can be exactly recovered with high probability via solving the convex optimization problem
\begin{equation}\label{prob:cpca}
  \min_{L,S} \left\{\|L\|_* + \lambda \|S\|_1: \; \hbox{s.t. } {\cal A}(L+S) = b\right\},
\end{equation}
as long as the range space of the adjoint operator ${\cal A}^*$ is randomly distributed according to the Haar measure and its dimension $q$ is in the order $O\left((\rho mn + mr)  \log^2m\right)$. Here $\lambda = 1/\sqrt{m}$ is a constant, $\|L\|_*$ and $\|S\|_1$ denote the nuclear norm of $L$ (sum of all singular values) and the $\ell_1$ norm of $S$ (sum of absolute values of all components), respectively. Note that to determine a rank $r$ matrix, it is sufficient to specify $(m+n-r)r$ elements. Let the number of nonzeros of $S_0$ be denoted by $\mathrm{nnz}(S_0)$. Without considering the distribution of the support of $S_0$, we define the \emph{degree of freedom} of the pair   $(L_0,S_0)$ by
\begin{equation}\label{def:dof}
\mathrm{dof} := (m+n-r)r + \mathrm{nnz}(S_0).
\end{equation}

The augmented Lagrangian function of \eqref{prob:cpca} is given by
\begin{equation*}
  {\cal \bar L}(L, S, p) := \|L\|_* + \lambda \|S\|_1 - \langle p, {\cal A}(L+S) - b\rangle + \frac{\beta}{2} \|{\cal A}(L+S) - b\|^2.
\end{equation*}
One can see that the minimization of ${\cal \bar L}$ with respect to either $L$ or $S$, with the other two variables being fixed, does not have closed  form  solution. To avoid inner loop for iteratively solving  ADMM-subproblems, the linearized ADMM framework \eqref{LADM-xpy} and its inertial version \eqref{iLADM-xpy} can obviously be applied. Note that it is necessary to linearize both ADMM-subproblems in order to avoid inner loops.
Though the iterative formulas of LADMM and inertial LADMM for solving \eqref{prob:cpca} can be derived very easily based on \eqref{LADM-xpy} and \eqref{iLADM-xpy}, we elaborate them below for clearness and subsequent references.
Let $(L^k,S^k,p^k)$ be given.
The LADMM framework \eqref{LADM-xpy} for solving \eqref{prob:cpca} appears as
\begin{subequations}
  \label{LADM-LyS}
\begin{eqnarray}
\label{LADM-LyS-U}
  U^k     &=& {\cal A}^*({\cal A} (L^k+ S^k)-b),\\
\label{LADM-LyS-L}
  L^{k+1} &=& \arg\min_{L} \|L\|_* - \langle p^k, {\cal A}(L) \rangle + \frac{\beta}{2\tau}\|L - (L^k -\tau U^k) \|_F^2, \\
\label{LADM-LyS-y}
  p^{k+1} &=& p^k - \beta ({\cal A} (L^{k+1}+S^k)-b), \\
\label{LADM-LyS-V}
  V^k &=& {\cal A}^*({\cal A} (L^{k+1}+S^k)-b),\\
\label{LADM-LyS-S}
  S^{k+1} &=& \arg\min_{S} \lambda\|S\|_1 - \langle p^{k+1}, {\cal A}(S)\rangle + \frac{\beta}{2\eta}\|S - (S^k -\eta V^k) \|_F^2.
\end{eqnarray}
\end{subequations}
The inertial LADMM framework \eqref{iLADM-xpy} for solving \eqref{prob:cpca} appears as
\begin{subequations}
  \label{iLADM-LyS}
\begin{eqnarray}
\label{iLADM-LyS-bar}
(\bar L^k,\bar S^k,\bar p^k) &=& (L^k, S^k, p^k) + \alpha_k (L^k - L^{k-1}, S^k - S^{k-1}, p^k - p^{k-1}), \\
\label{iLADM-LyS-U}
  U^k     &=& {\cal A}^*({\cal A} (\bar L^k+ \bar S^k)-b),\\
\label{iLADM-LyS-L}
  L^{k+1} &=& \arg\min_{L} \|L\|_* - \langle \bar p^k, {\cal A}(L) \rangle + \frac{\beta}{2\tau}\|L - (\bar L^k -\tau U^k) \|_F^2, \\
\label{iLADM-LyS-y}
  p^{k+1} &=& \bar p^k - \beta ({\cal A} (L^{k+1}+\bar S^k)-b), \\
\label{iLADM-LyS-V}
  V^k &=& {\cal A}^*({\cal A} (L^{k+1}+\bar S^k)-b),\\
\label{iLADM-LyS-S}
  S^{k+1} &=& \arg\min_{S} \lambda\|S\|_1 - \langle p^{k+1}, {\cal A}(S)\rangle + \frac{\beta}{2\eta}\|S - (\bar S^k -\eta V^k) \|_F^2.
\end{eqnarray}
\end{subequations}
Note that the subproblems \eqref{LADM-LyS-L} (or \eqref{iLADM-LyS-L}) and \eqref{LADM-LyS-S} (or \eqref{iLADM-LyS-S}) have closed form solutions given, respectively,
by the shrinkage operators of matrix nuclear norm and vector $\ell_1$ norm, see, e.g., \cite{MGC11,YZ11}. The main computational cost per iteration of both algorithms
is one singular value decomposition (SVD) required in solving  the $L$-subproblem.

\subsection{Generating experimental data}
In our experiments, we set $m=n$ and tested different ranks of $L_0$ (denoted by $r$), sparsity levels of $S_0$ (i.e., $\text{nnz}(S_0)/(mn)$) and sample ratios (i.e., $q/(mn)$). The low-rank matrix $L_0$ was generated by $\text{randn}(m,r)*\text{randn}(r,n)$ in MATLAB. The support of $S_0$ is randomly determined by uniform distribution, while the values of its nonzeros are uniformly distributed in $[-10,10]$. Such type of synthetic data are roughly those tested in \cite{WGMM13}.
As for the linear operator ${\cal A}$, we tested three types of linear operators, i.e., two-dimensional partial DCT (discrete cosine transform), FFT (fast Fourier transform) and WHT (Walsh-Hadamard transform). The rows of these transforms are selected uniformly at random.

\subsection{Parameters, stopping criterion and initialization}
The model parameter $\lambda$ was set to $1/\sqrt{m}$ in our experiments, which is determined based on the exact recoverability theory in \cite{WGMM13}.
As for the other parameters ($\beta$, $\tau$ and $\eta$) common to LADMM and iLADMM, we used the same set of values and adaptive rules in all the tests. Now we elaborate how the parameters are chosen.
Since ${\cal A}$ contains rows of orthonormal transforms, it holds that ${\cal A}{\cal A}^* = {\cal I}$, the identity operator. Therefore,  it holds that $\rho({\cal A}^*{\cal A})=1$. We set $\tau = \eta = 0.99$, which satisfies the convergence requirement specified in Theorems \ref{Theorem-convergence-LADM} and  \ref{Theorem-iLADM-xpy}. The penalty parameter $\beta$ was initialized at $0.1q/\|b\|_1$ and was tuned at the beginning stage of the algorithm.
Specifically,  we tuned $\beta$  within the first $30$ iterations according to the following rule:
\[
\beta_{k+1} =
         \left\{
          \begin{array}{ll}
            \max(0.5\beta_k,10^{-3}), & \hbox{if $r_k < 0.1$;} \\
            \min(2\beta_k, 10^2), & \hbox{if $r_k > 5$;} \\
            \beta_k, & \hbox{otherwise,}
          \end{array}
        \right.
\text{~~where~~} r_k := \frac{\beta_k\|{\cal A}(L^k+S^k)-b\|^2}{2s_k (\|L^k\|_* + \lambda \|S^k\|_1)}.
\]
Here $s_k$ is a parameter attached to the objective function $\|L\|_* + \lambda\|S\|_1$ and was chosen adaptively so that the quadratic term $\frac{\beta}{2}\|{\cal A}(L+S)-b\|^2$ and the objective term $\|L\|_* + \lambda \|S\|_1$ remain roughly in the same order.
Note that the choice of $\beta$ does not have much theory and is usually determined via numerical experiments, see, e.g., \cite{YY13} for the influence of different $\beta$'s in linearized ADMM for matrix completion problem.
The extrapolation parameter $\alpha_k$ for iLADMM was set to $0.28$ and held constant in all our experiments. Note that this value of $\alpha_k$ is determined based on  experiments and may be far from optimal. How to select $\alpha_k$ adaptively to achieve stable and faster convergence remains a research issue. Here our main goal is to illustrate the effect of the extrapolation steps.
We also present some numerical results to compare the performance of iLADMM with different constant strategies for $\alpha_k$.

It is easy to see from \eqref{LADM-mVI-k+1} that if two consecutive iterates generated by proximal point method  are identical then a solution is already obtained.
Since LADMM is an application of a general PPA,  we terminated it by the following rule
\begin{equation}\label{stop-rule}
  \frac{\|(L^{k+1},S^{k+1},p^{k+1}) - (L^k,S^k,p^k)\|}{1 + \|(L^k,S^k,p^k)\|} < \varepsilon,
\end{equation}
where $\varepsilon > 0$ is a tolerance parameter. Here $\|(L,S,p)\| := \sqrt{\|L\|_F^2 + \|S\|_F^2 + \|p\|^2}$.
Since iLADMM generates the new point $(L^{k+1},S^{k+1},p^{k+1})$ by applying proximal point method to  $ (\bar L^k,\bar S^k,\bar p^k)$,
we used the same stopping rule as \eqref{stop-rule} except that $(L^k,S^k,p^k)$ is replaced by $ (\bar L^k,\bar S^k,\bar p^k)$. That is
\begin{equation}\label{stop-rule-iLADMM}
  \frac{\|(L^{k+1},S^{k+1},p^{k+1}) - (\bar L^k,\bar S^k,\bar p^k)\|}{1 + \|(\bar L^k,\bar S^k,\bar p^k)\|} < \varepsilon.
\end{equation}
%
In our experiments, we initialize all variables $L$, $S$ and $p$ at zeros. 

\subsection{Experimental results}
Recall that  the matrix size is $m\times n$, the number of measurements is $q$, the rank of $L_0$ is $r$, and the degree of freedom of the pair $(L_0,S_0)$ is defined in \eqref{def:dof}.
In our experiments, we tested $m=n=1024$.
Let $k$ be the number of nozeros of $S_0$.
We tested four different ranks for $L_0$, three levels of sparsity for $S_0$ and four levels of sample ratios. Specifically, in our experiments
we tested  $r\in \{5,10,15,20\}$, $k/m^2 \in \{0.01, 0.05, 0.10\}$ and $q/m^2 \in \{0.4,0.6,0.8\}$.

Let $(L,S)$ be the recovered solution. For each setting, we report the relative errors of $L$ and $S$ to the true low-rank and sparse matrices $L_0$ and $S_0$, i.e.,
$\|L-L_0\|_F/\|L_0\|_F$ and $\|S-S_0\|_F/\|S_0\|_F$, and the number of iterations to meet the condition \eqref{stop-rule} or \eqref{stop-rule-iLADMM}, which are denoted by iter1 and iter2 for LADMM and iLADMM, respectively. We terminated both algorithms if the number of iterations reached 1000 but the stopping rule \eqref{stop-rule} or \eqref{stop-rule-iLADMM} still did not hold.
For each problem scenario, we run 10 random trials for both algorithms and report the averaged results. Detailed experimental results for $\varepsilon=10^{-5}$ and $r=5,10,15$ and $20$ are given in Tables \ref{table-r=5}-\ref{table-r=20}, respectively. In each table, dashed line  ``---" represents that the maximum iteration number was reached.

\begin{table}\caption{Results of $\text{rank}(L_0) = 5$: $\varepsilon=10^{-5}$, average results of 10 random trials.}
\begin{center}
\begin{footnotesize}
\input{results/TableHead.tex}
\input{results/18-Jul-2014/r_5.dat}
\end{tabular}\label{table-r=5}
\end{footnotesize}
\end{center}
\end{table}

 \begin{table}\caption{Results of $\text{rank}(L_0) = 10$: $\varepsilon=10^{-5}$, average results of 10 random trials.}
\begin{center}
\begin{footnotesize}

\input{results/TableHead.tex}
\input{results/18-Jul-2014/r_10.dat}
\end{tabular}\label{table-r=10}
\end{footnotesize}
\end{center}
\end{table}

 \begin{table}\caption{Results of $\text{rank}(L_0) = 15$: $\varepsilon=10^{-5}$, average results of 10 random trials.}
\begin{center}
\begin{footnotesize}

\input{results/TableHead.tex}
\input{results/18-Jul-2014/r_15.dat}
\end{tabular}\label{table-r=15}
\end{footnotesize}
\end{center}
\end{table}

\begin{table}\caption{Results of $\text{rank}(L_0) = 20$: $\varepsilon=10^{-5}$, average results of 10 random trials.}
\begin{center}
\begin{footnotesize}

\input{results/TableHead.tex}
\input{results/18-Jul-2014/r_20.dat}
\end{tabular}\label{table-r=20}
\end{footnotesize}
\end{center}
\end{table}

It can be seen from Tables \ref{table-r=5}-\ref{table-r=20} that iLADMM is generally faster than LADMM to obtain solutions satisfying the aforementioned conditions.  Specifically,  within our setting  the numbers of iterations consumed by iLADMM range, roughly, from $60\%$--$80\%$ of those consumed by LADMM. If we take into account all the tests (except those cases where either LADMM or iLADMM failed to terminate within 1000 iterations, e.g., $(r,k/m^2,q/m^2) = (5,0.1,40\%)$ and ${\cal A}$ is partial DCT), the overall average number of iterations used by iLADMM is about $74\%$ of that used by LADMM. Note that in some cases iLADMM  obtained satisfactory results within the number of allowed iterations (1000 in our setting), while LADMM did not. For example, $(r,k/m^2,q/m^2) = (5,0.1,40\%)$ and ${\cal A}$ is partial DCT or partial WHT. In most  cases, the recovered matrices $L$ and $S$ are close to the true low-rank and sparse components $L_0$ and $S_0$, respectively. The relative errors are usually in the order $10^{-5}$---$10^{-6}$. For some cases, the recovered solutions are not of high quality (relative errors are large), which is mainly because the number of samples are small relative to the degree of freedom of $(L_0,S_0)$. This can be seen from the values of $q/\text{dof}$ listed in the tables. Roughly speaking, the recovered solutions are satisfactory (say, relative errors are less than $10^{-3}$) provided that $q/\text{dof}$ is no less than $3.5$.

We note that the per iteration cost of both LADMM and iLADMM for the compressive principal pursuit model \eqref{prob:cpca} is dominated by one SVD and thus is roughly identical. The extra cost of the extrapolation inertial step in \eqref{iLADM-LyS-bar}  is negligible compared to the computational load of SVD. This is the main reason that we only reported the number of iterations but not CPU time consumed by both algorithms. The inertial technique actually accelerates the original algorithm to a large extent but without increasing the total computational cost.

\begin{figure}[htbp]
\centering{
\includegraphics[trim = 0 200 50  160, scale = 0.27]{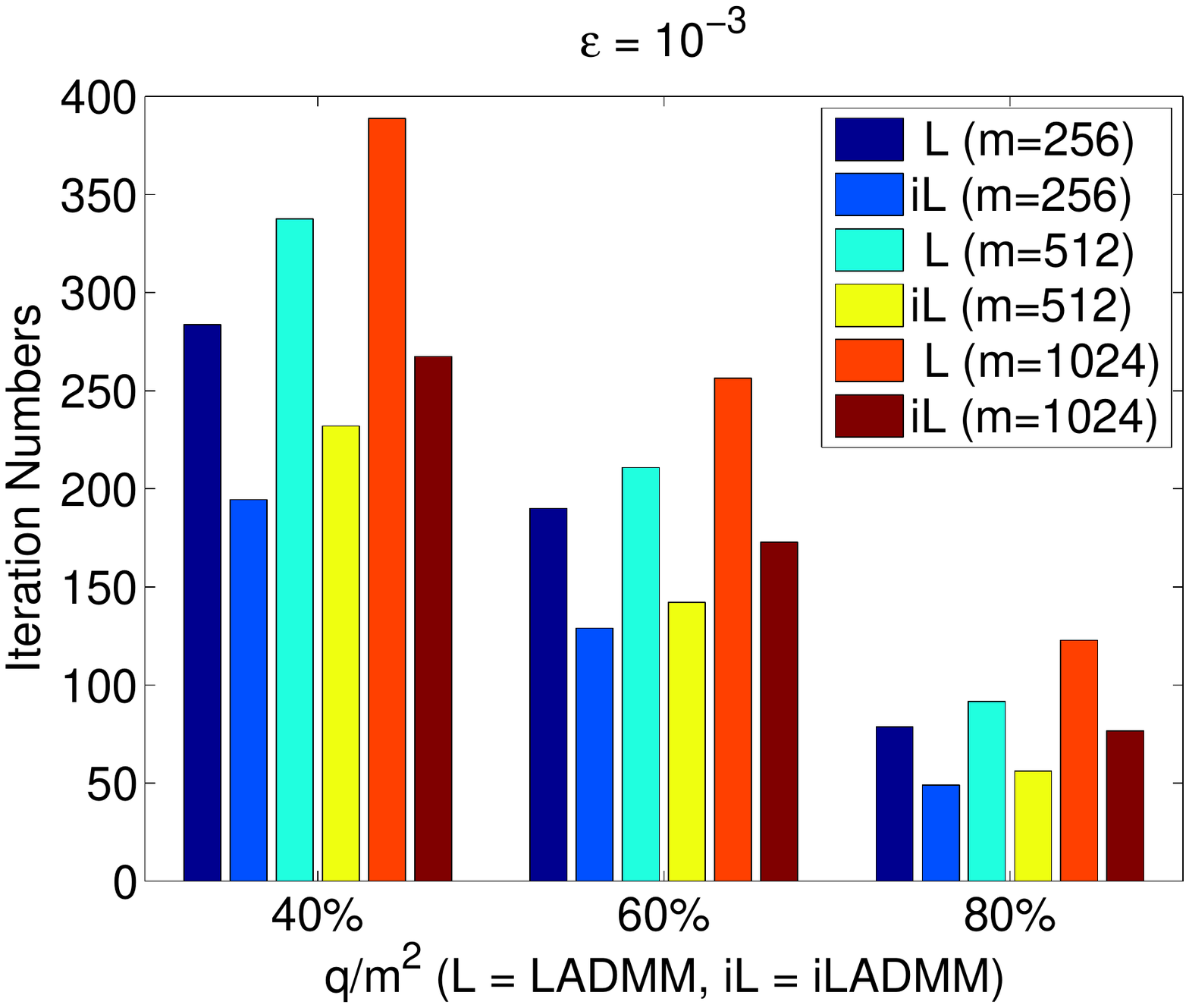}
\includegraphics[trim = 0 200 50  160, scale = 0.27]{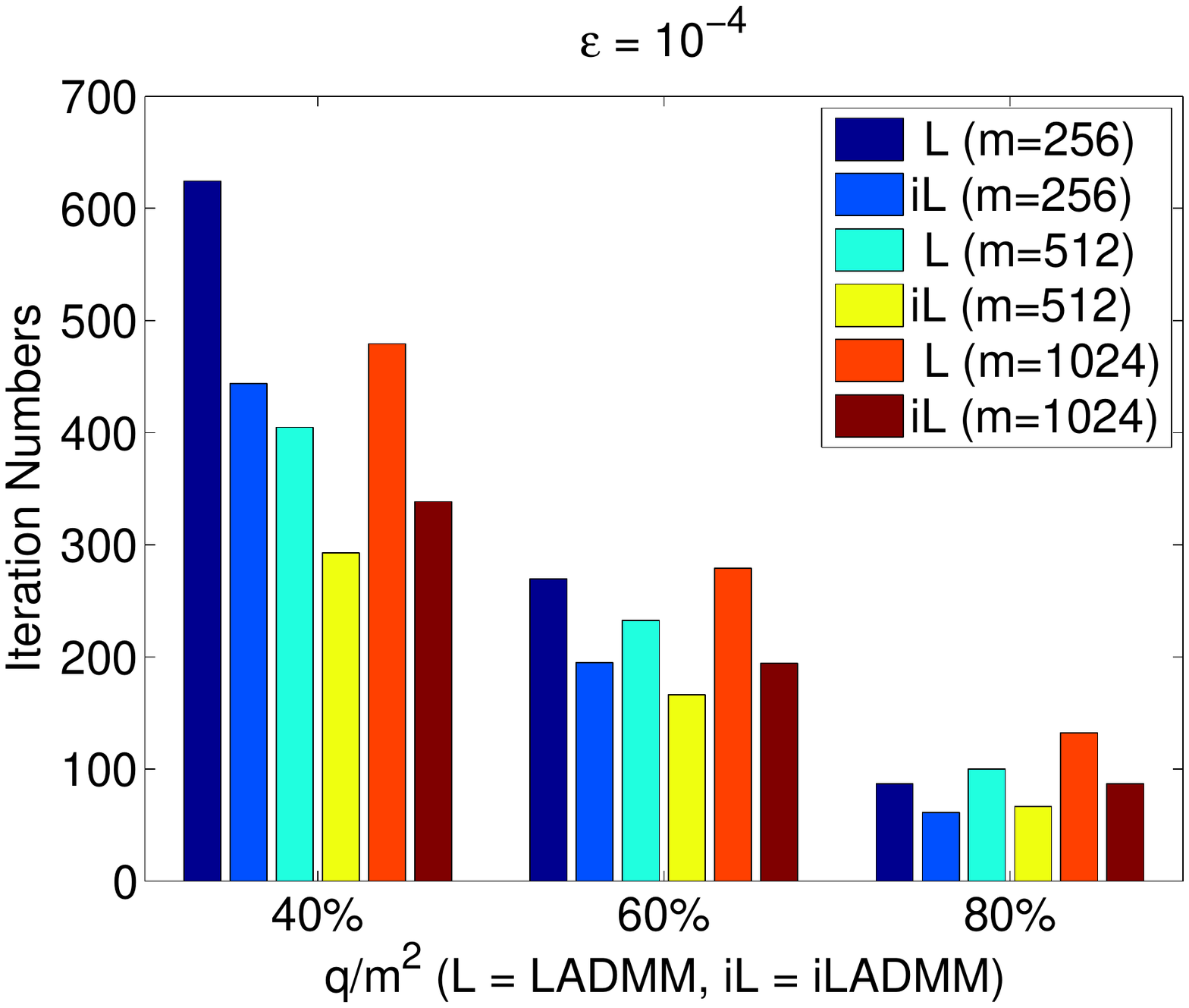}
\includegraphics[trim = 0 200 50  160, scale = 0.27]{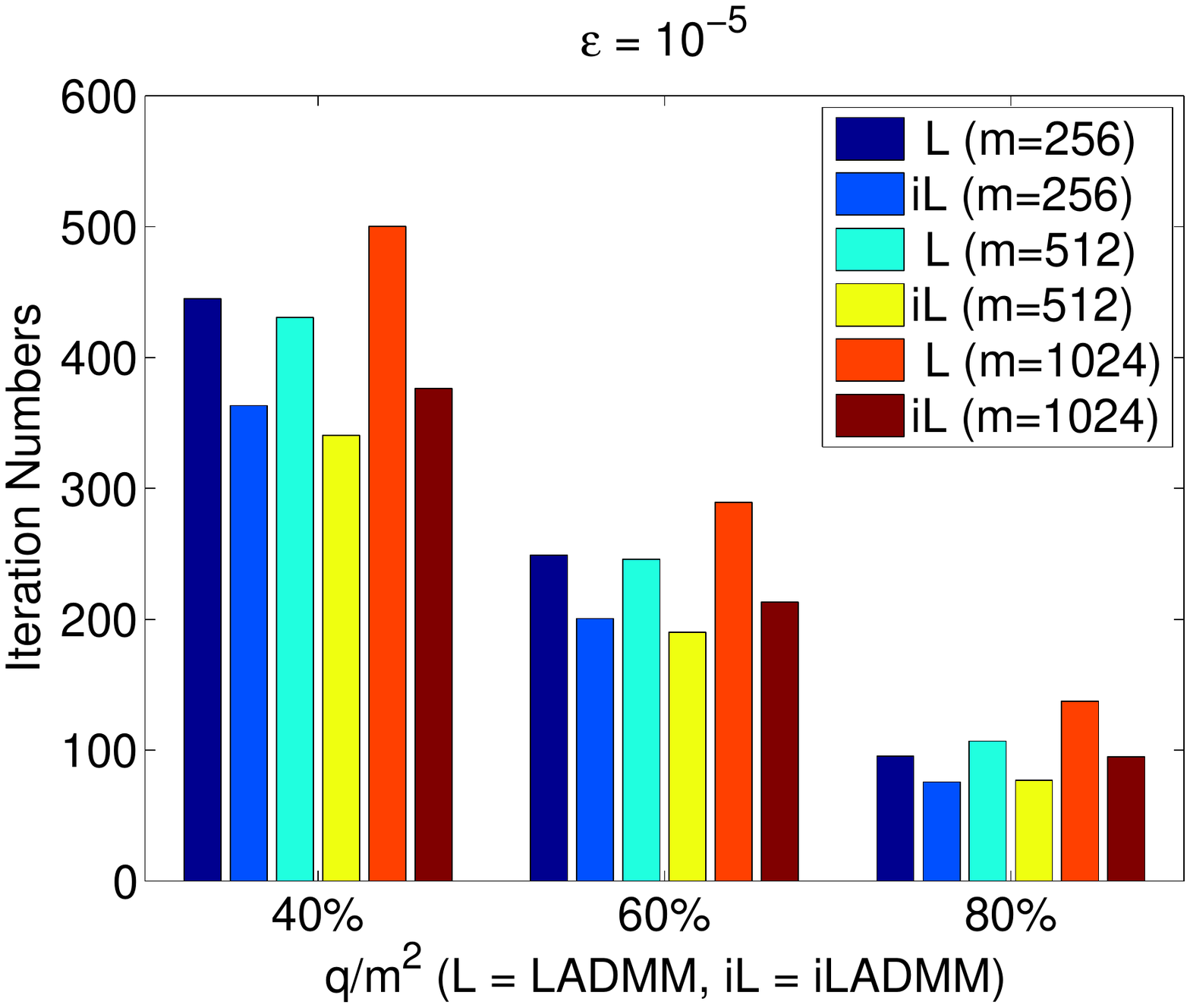}
}\caption{Comparison results on different matrix sizes  and stopping tolerance:  Average results of 10 random trials ($m=n=256,512,1024$, and from left to right $\varepsilon = 10^{-3}, 10^{-4}, 10^{-5}$, respectively). }
\label{Fig-m3}
\end{figure}

To better understand the behavior of iLADMM relative to LADMM, we also tested different matrix sizes ($m=n=256, 512$ and $1024$) with different levels of stopping tolerance ($\varepsilon=10^{-3},10^{-4}$ and $10^{-5}$ in \eqref{stop-rule}). For each case, we tested $r\in\{5,10,15,20\}$ and $k/m^2\in\{0.01,0.05,0.10\}$ for a fixed $q$ such that $q/m^2 \in \{0.4,0.6,0.8\}$.
For each $q$, we accumulated the iteration numbers for different $(r,k)$ and the three types of linear operators and took an average finally. The results are summarized in Figure \ref{Fig-m3}. Again, these results are average of 10 random trials for each case. From the results we can see that iLADMM is faster and terminates earlier than LADMM with different levels of stopping tolerance. Roughly speaking, iLADMM reduced the cost of LADMM by about 30\%.

\begin{figure}[htbp]
\centering{
\includegraphics[trim = 0 200 50  160, scale = 0.27]{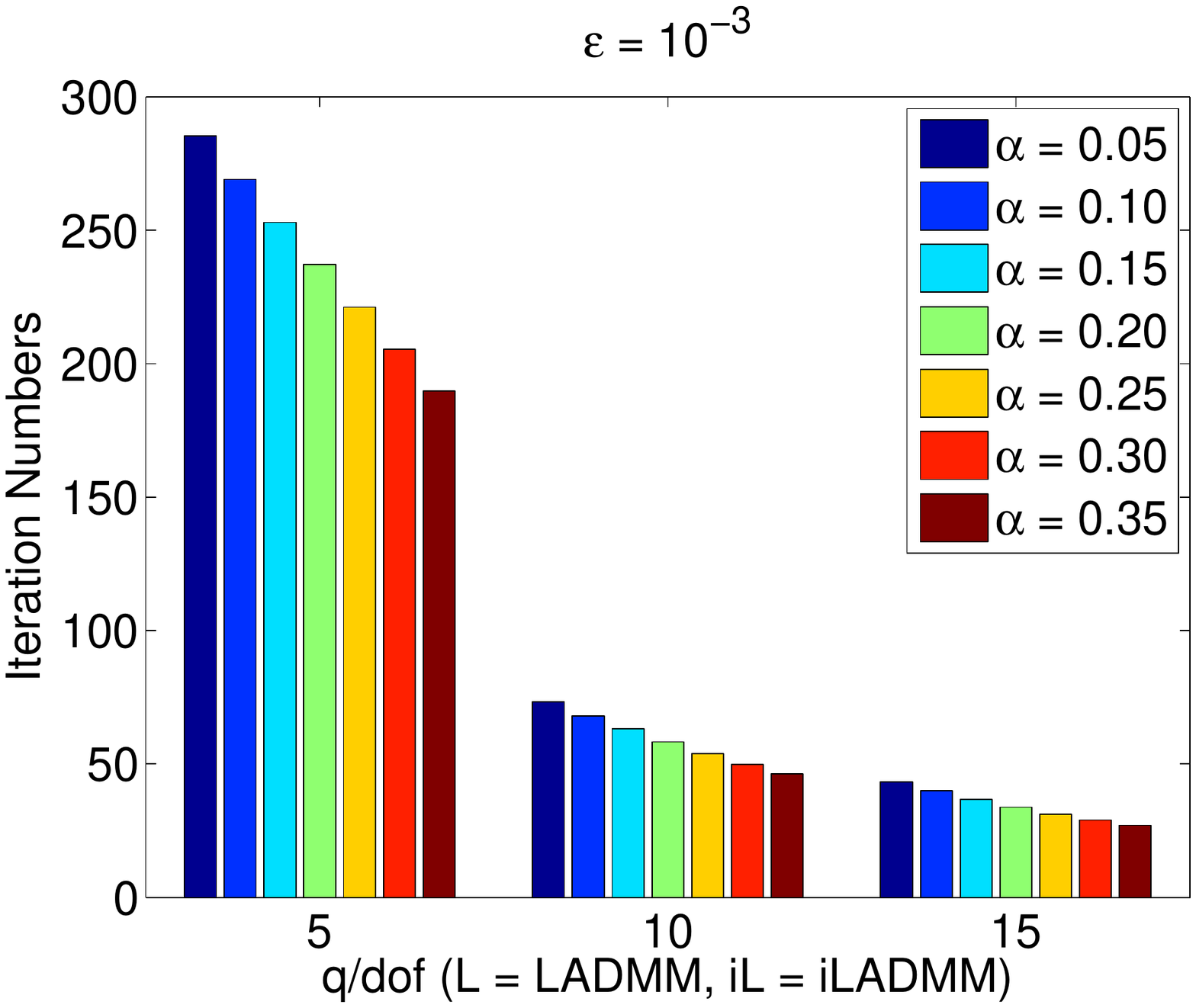}
\includegraphics[trim = 0 200 50  160, scale = 0.27]{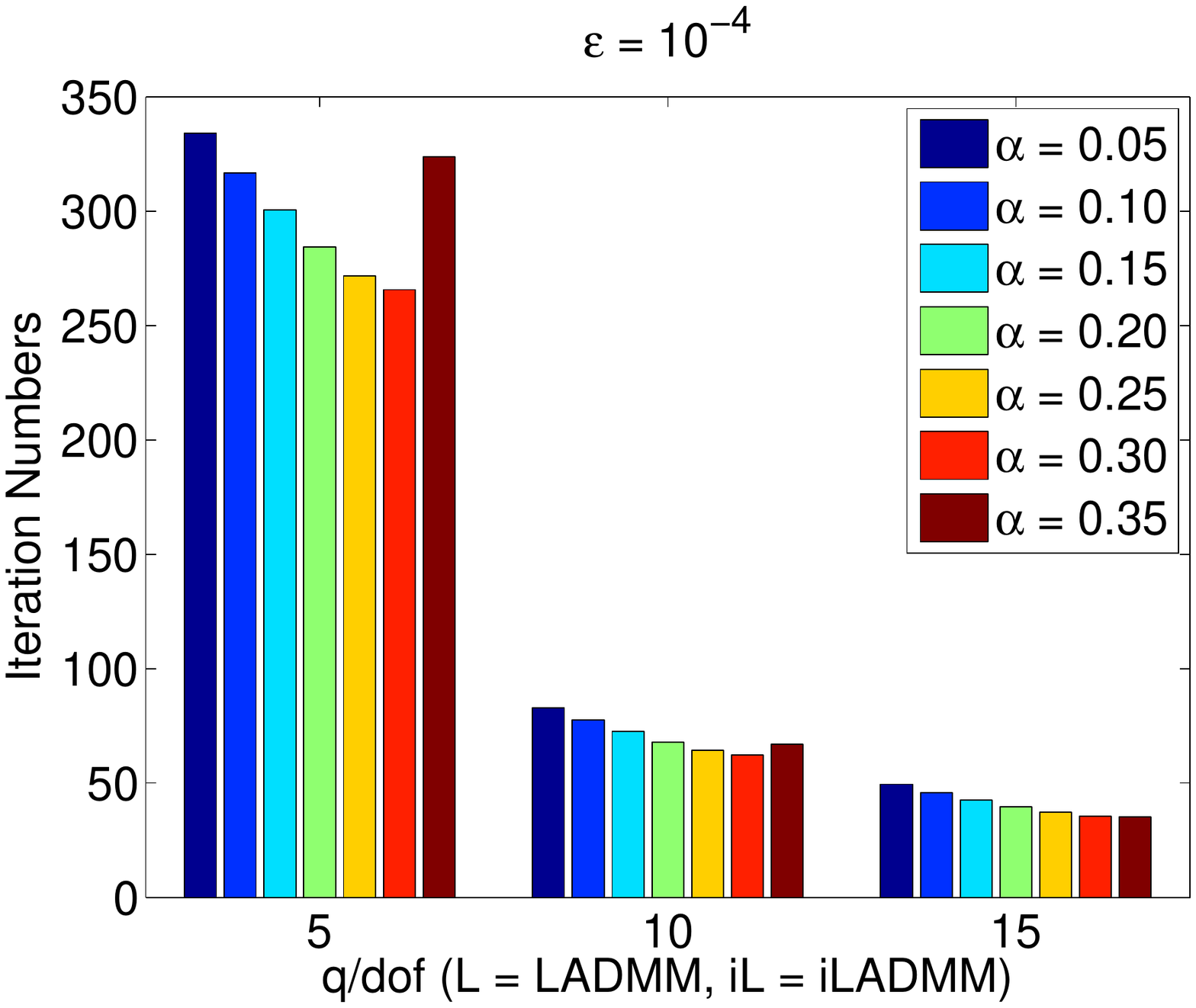}
\includegraphics[trim = 0 200 50  160, scale = 0.27]{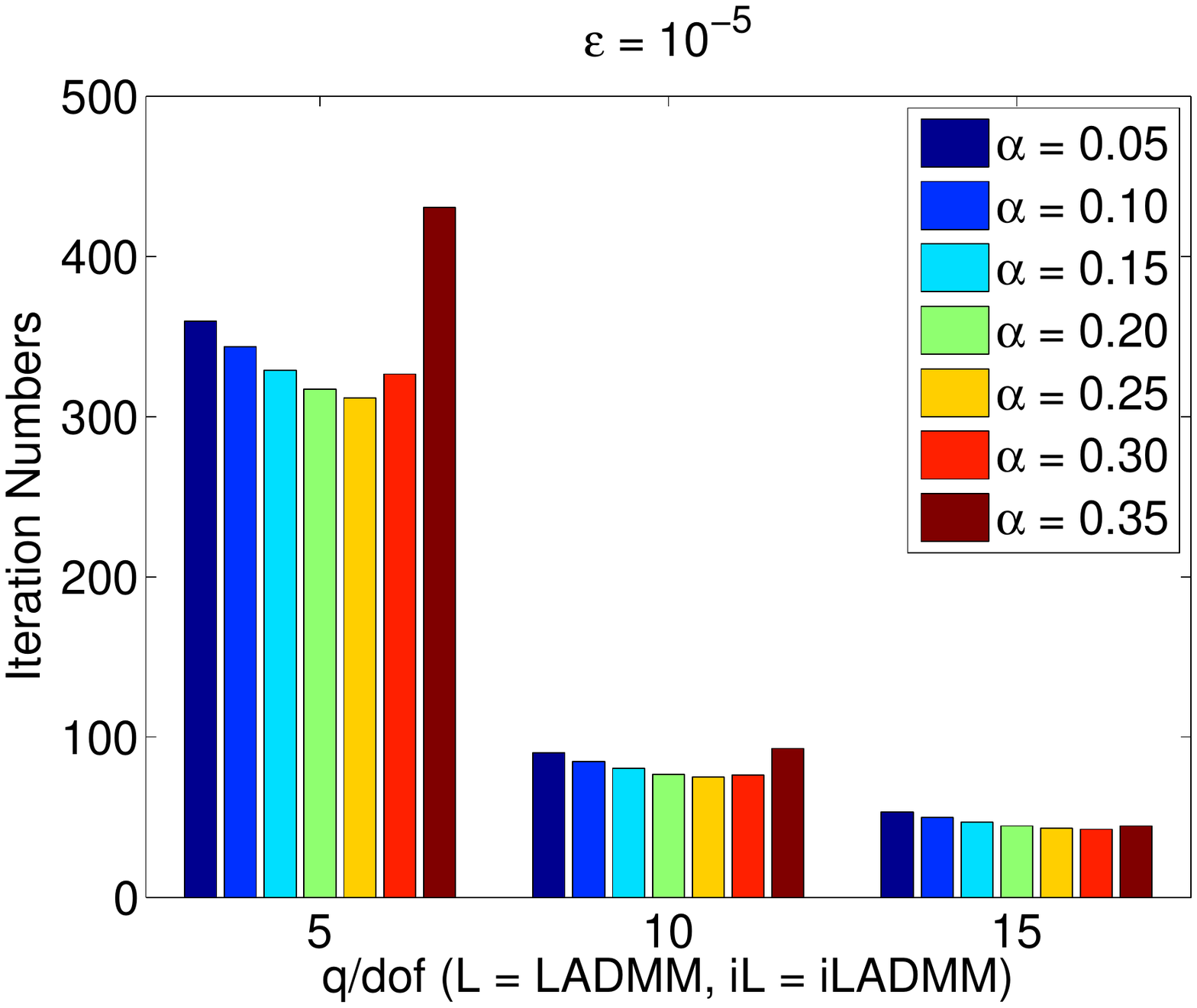}
}\caption{Comparison results on different $\alpha_k\equiv\alpha$ and stopping tolerance: Average results of 10 random trials ($m=n=512$, $\alpha$ ranges from $0.05$ to $0.35$,
and from left to right  $\varepsilon = 10^{-3}, 10^{-4}, 10^{-5}$, respectively).}
\label{Fig-alp6}
\end{figure}

We also run iLADMM with various constant strategies for $\alpha_k$. In particular, we set $m=n=512$ and tested different values of  $q$ such that $q/\text{dof}\in \{5,10,15\}$. For each case, we varied $r \in \{5,10,15,20\}$ and $k/m^2 \in \{0.01,0.05,0.10\}$ for the three types of aforementioned measurement matrices. We accumulated the number of iterations and took an average finally. The detailed average results of 10 random trials for $\alpha_k\equiv \alpha$  from $0.05$ to $0.35$ are given in Figure \ref{Fig-alp6}.

From the results in Figure \ref{Fig-alp6} we see that, for the tested 7 values of $\alpha$, iLADMM is slightly faster if $\alpha$ is larger, provided that $\alpha$ does not exceed $0.3$.  %
We have also observed that for $\alpha > 0.3$ iLADMM either slows down or performs not very stable, especially when $q/\text{dof}$ is small. This is the main reason that we set $\alpha_k$ a constant value that is near $0.3$ but not larger.

\section{Concluding remarks}
In this paper, we proposed and analyzed a general inertial proximal point method within the setting of mixed VI problem \eqref{mVI}. The proposed method adopts a weighting matrix and allows more flexibility.
Our convergence results require weaker conditions in the sense that the weighting matrix $G$ does not necessarily be positive definite, as long as the function $F$ is $H$-monotone and $G$ is positive definite in the null space of $H$. The convergence analysis can be easily adapted to the monotone inclusion problem \eqref{prob:0inTw}.
We also  showed  that the linearized ADMM for linearly constrained separable convex optimization problem is a proximal point method applied to the primal-dual optimality conditions, as long as the parameters are reasonably small. As byproducts of this finding, we  established  with standard analytic techniques for proximal point method the global convergence and convergence rate results of LADMM.
This proximal reformulation also allows us to propose an inertial version of LADMM, whose convergence is guaranteed under suitable conditions. Our preliminary implementation of the algorithm and extensive experimental results on compressive principal component pursuit problem have shown that the inertial LADMM is generally faster than the original LADMM. Though in a sense the acceleration is not very significant, we note that the inertial LADMM does not require any additional and unnegligible computational cost either.

Throughout our experiments  the extrapolation steplength $\alpha_k$ held constant. How to select $\alpha_k$ adaptively based on the current information such that the overall algorithm performs more efficiently and stable is a practically very important question and deserves further investigation. Another theoretical issue is to investigate  worst-case complexity analysis for general inertial type algorithms.
In fact, complexity results of inertial type algorithms for minimizing closed proper convex functions already exist in the literature. The pioneering work in this direction is due to Nesterov \cite{Nes83}, where the algorithm can also be viewed in the perspective of inertial algorithms. Refined analyses for more general problems can be found in \cite{BT09,Gul92}.
Let $f: \Re^n\rightarrow \Re$ be a closed proper convex function and be bounded below.
Based on   \cite{Nes83,BT09,Gul92}, the following algorithm can be studied. Let $x^0 \in\Re^n$ be given. Set $x_0 =  x^{-1}$, $t_0=1$ and $k=0$. For $k\geq 0$  the algorithm iterates as
\begin{subequations}\label{Alg-optimal-iPPA}
  \begin{eqnarray}
  \label{Alg-optimal-t}
    t_{k+1} &=& \frac{1+\sqrt{1+4t_k^2}}{2},  \\
  \label{Alg-optimal-bar-w}
    \bar w^k &=& w^k + \frac{t_k-1}{t_{k+1}} (w^k - w^{k-1}), \\
  \label{Alg-optimal-w}
    w^{k+1} &=& \arg\min_w f(w) + \frac{1}{2\lambda_k}\|w - \bar w^k\|^2.
  \end{eqnarray}
\end{subequations}
Using analyses similar to those in \cite{Nes83, Gul92, BT09}, one can show that the sequence $\{w^k\}_{k=0}^\infty$ satisfies
\[f(w^k) - \min_{w\in\Re^n} f(w) = O(1/k^2).\]
Algorithm \eqref{Alg-optimal-iPPA} is nothing but an inertial PPA with steplength $\alpha_k = \frac{t_k-1}{t_{k+1}}$. It is interesting to note that $\alpha_k$ is monotonically increasing as $k\rightarrow \infty$ and converges to $1$, which is much larger than the upper bound condition $\alpha<1/3$ required in Theorem \ref{Theorem2}.
Also note that the convergence for \eqref{Alg-optimal-iPPA} is measured by the objective residue. Without further assumptions on $f$, it seems difficult to establish convergence of the sequence $\{w^k\}_{k=0}^\infty$, see, e.g., \cite{Gul92}.
In comparison, our results impose smaller upper bound on $\alpha_k$ but guarantee the convergence of the sequence of iterates $\{w^k\}_{k=0}^{\infty}$.
Even though, there seems to be certain gap between the classical results \cite{Nes83,BT09,Gul92} for minimizing closed proper convex functions and the results presented in the present paper. Further research in this direction is interesting.

\def\cprime{$'$}

\begin{appendix}

\section{Proof of Theorem~\ref{Theorem-convergence-LADM}}
First, we sketch the proof of convergence of the sequence $\{w^k\}$ to a solution of \eqref{mVI}.
Clearly, the matrix $G$ defined in \eqref{G-LADM} is symmetric and positive definite under our assumptions that $\tau < 1/\rho(A^TA)$ and $\eta < 1/\rho(B^TB)$.
Let $w^*\in \Omega^*$ be arbitrary fixed. It follows from setting $w=w^*$ in \eqref{LADM-mVI-k+1} that
\begin{eqnarray*}
  \langle w^{k+1} - w^*, G(w^k - w^{k+1})\rangle &\geq& \theta(w^{k+1}) - \theta(w^*) + \langle w^{k+1}-w^*, F(w^{k+1})\rangle \\
  &\geq& \theta(w^{k+1}) - \theta(w^*) + \langle w^{k+1}-w^*, F(w^*)\rangle \\
  &\geq& 0,
\end{eqnarray*}
where the second ``$\geq$" follows form the monotonicity of $F$. Therefore, we obtain
\begin{eqnarray}
  \label{PPA-decrease} \nonumber
  \|w^{k+1}-w^*\|_G^2 &=& \|w^k-w^*\|_G^2 - \|w^k-w^{k+1}\|_G^2 - 2 \langle w^{k+1}-w^*, G(w^k-w^{k+1})\rangle \\
  &\leq& \|w^k-w^*\|_G^2 - \|w^k-w^{k+1}\|_G^2.
\end{eqnarray}
Since $G$ is positive definite, this implies that measured by $G$-norm the sequence $\{w^k\}$ is strictly contractive with respect to $\Omega^*$ unless $w^k=w^{k+1}$ in which case $w^k$ is already a solution. The convergence of $\{ w^k \}$ to some solution $w^\star\in \Omega^*$ follows directly from standard analyses for PPA and the key inequality \eqref{PPA-decrease}. We omit the details.

 Second, we prove \eqref{ergodic-rate}.
Let $w^{i+1} \in \Omega$ be generated via \eqref{LADM-xpy}. It follows from the monotonicity of $F$ and \eqref{LADM-mVI-k+1} that, for any $w\in \Omega$, there holds
\begin{eqnarray*}
  \theta(w) - \theta(w^{i+1}) + (w - w^{i+1})^T F(w)    &\geq&  \theta(w) - \theta(w^{i+1}) + (w - w^{i+1})^T F(w^{i+1}) \\
  &\geq& (w - w^{i+1})^T G ( w^i - w^{i+1}).
\end{eqnarray*}
By noting the relation
 $2(w - w^{i+1})^T G ( w^i - w^{i+1}) \geq  \|w - w^{i+1}\|_G^2 -\|w - w^i\|_G^2 $,
 we obtain
\[
 \theta(w) - \theta(w^{i+1}) + (w - w^{i+1})^T F(w)    \geq
 \frac{1}{2}\left( \|w - w^{i+1}\|_G^2  -\|w - w^i\|_G^2  \right), \forall w\in \Omega.
 \]
Take sum over $i = 0,1,\ldots,k$ and  divide both sides by $(k+1)$, we get
 \begin{equation}
   \label{O(1/k)-final}
    \theta(w) - \frac{1}{k+1}\sum_{i=0}^k\theta(w^{i+1}) + \left(w - \frac{1}{k+1}\sum_{i=0}^k w^{i+1}\right)^T F(w)
 \geq
 - \frac{\|w - w^0\|_G^2}{2(k+1)}, \;
 \forall w\in \Omega.
  \end{equation}
The conclusion \eqref{ergodic-rate} follows directly from \eqref{O(1/k)-final} by noting the definition of $\bar w^k$ and the fact that
  \[
  \frac{1}{k+1}\sum_{i=0}^k \theta(w^{i+1}) \geq \theta\left(\frac{1}{k+1}\sum_{i=0}^k w^{i+1}\right) =  \theta(\bar w^k).
 \]
The equivalence of \eqref{saddle-conditions-approximate} and \eqref{ergodic-rate} can be verified directly from the notation defined in \eqref{def-w-F-theta} and the definition of ${\cal L}$ in \eqref{def:L}.

 Finally, we prove \eqref{non-ergodic-rate} and \eqref{non-ergodic-rate-o}. Since \eqref{LADM-mVI-k+1} holds for all $k$, it also holds for $k := k-1$, i.e.,
  \begin{equation}\label{LADM-mVI-k}
    \theta(w) - \theta(w^{k})   + \langle w - w^{k},   F(w^{k})   + G (w^{k} -  w^{k-1})\rangle \geq 0, \; \forall w\in\Omega.
  \end{equation}
 By setting $w=w^k$ and $w=w^{k+1}$ in \eqref{LADM-mVI-k+1} and \eqref{LADM-mVI-k}, respectively, and taking an addition, we obtain
  \[
        \langle G(w^{k+1} - w^{k}),    (w^{k} -  w^{k-1}) -  (w^{k+1} -  w^k)\rangle  \geq  \langle w^{k+1} - w^{k}, F(w^{k+1}) - F(w^{k})\rangle \geq 0.
  \]
In addition, by taking into account the fact that
  \begin{eqnarray*}
   \|w^k - w^{k-1}\|_G^2-\|w^{k+1} - w^k\|_G^2  \geq  2\langle G(w^{k+1}-w^k), (w^k - w^{k-1})-(w^{k+1}-w^k)\rangle,
  \end{eqnarray*}
  we obtain $\|w^k - w^{k-1}\|_G \geq \|w^{k+1} - w^k\|_G$, i.e.,  $\|w^k - w^{k-1}\|_G$ is monotonically nonincreasing with respect to $k$.
By further considering  \eqref{PPA-decrease}, we obtain
  \[
  k\|w^{k} - w^{k-1}\|_G^2 \le \sum_{i=0}^{k-1}\| w^{i+1} - w^{i}\|_G^2 \leq  \sum_{i=0}^{k-1} \left(\|w^i - w^*\|_G^2 -\|w^{i+1} - w^*\|_G^2\right) \leq \|w^0 - w^*\|_G^2,
  \]
which implies the relation \eqref{non-ergodic-rate}.
By using the trick introduced in \cite{DLY13,DY14},
we can derive the $o\left(1/k\right)$ result \eqref{non-ergodic-rate-o}. Specifically, we have
\begin{equation}\label{ch-add3}
  {k\over 2} \|w^k-w^{k-1}\|_G^2 \leq \sum_{i =\lfloor {k\over 2} \rfloor}^k \|w^{i} -w^{i-1}\|_G^2,
\end{equation}
where $\lfloor {k/2} \rfloor$ denotes the greatest integer no greater than $k/2$.
The result \eqref{non-ergodic-rate-o} follows by further considering
$\sum_{k=0}^{\infty}\| w^{k+1} - w^{k}\|_G^2 < \infty$ and thus the right-hand-side of \eqref{ch-add3}
converges to $0$ as $k\rightarrow\infty$.

\end{appendix}

\end{document}

%% file: results/TableHead.tex

\begin{tabular}{|c|c|c|c||c|c|r||c|c|r||c|}\hline
\mc{4}{|c||}{$m=n=1024$} & \mc{3}{c||}{LADMM} & \mc{3}{c||}{iLADMM} & \\ [1pt]\cline{1-10}
\mc{1}{|c|}{$r$} & \mc{1}{c|}{$k/m^2$ }  & \mc{1}{c|}{$(q/m^2, q/\text{dof})$}  & \mc{1}{c||}{${\cal A}$}
& \mc{1}{c|}{$\frac{\|L-L_0\|_F}{\|L_0\|_F}$}  & \mc{1}{c|}{$\frac{\|S-S_0\|_F}{\|S_0\|_F}$} & \mc{1}{c||}{iter1}
& \mc{1}{c|}{$\frac{\|L-L_0\|_F}{\|L_0\|_F}$}  & \mc{1}{c|}{$\frac{\|S-S_0\|_F}{\|S_0\|_F}$} & \mc{1}{c||}{iter2}
& $\frac{\text{iter2}}{\text{iter1}}$ \\  [2pt] \hline\hline

%% file: results/18-Jul-2014/r_5.dat
$5$   &   1\%   &   $(40\%,20.26)$  &  pdct & 2.50e-5 & 4.73e-5 & 269.6 & 1.47e-5 & 3.66e-5 & 205.4 & 0.76\\\cline{4-11}

   &      &     &  pfft & 8.71e-6 & 2.50e-5 & 220.4 & 1.36e-5 & 4.34e-5 & 170.8 & 0.77\\\cline{4-11}

   &      &     &  pwht & 3.13e-5 & 5.81e-5 & 263.9 & 1.81e-5 & 4.26e-5 & 201.9 & 0.76\\\cline{3-11}

   &      &   $(60\%,30.39)$  &  pdct & 1.23e-5 & 2.32e-5 & 189.1 & 1.54e-5 & 2.30e-5 & 136.9 & 0.72\\\cline{4-11}

   &      &     &  pfft & 9.32e-6 & 1.95e-5 & 148.2 & 1.32e-5 & 2.85e-5 & 105.4 & 0.71\\\cline{4-11}

   &      &     &  pwht & 9.53e-6 & 2.04e-5 & 159.6 & 8.54e-6 & 2.06e-5 & 117.9 & 0.74\\\cline{3-11}

   &      &   $(80\%,40.52)$  &  pdct & 3.02e-6 & 7.05e-6 & 90.9 & 3.86e-6 & 1.15e-5 & 61.3 & 0.67\\\cline{4-11}

   &      &     &  pfft & 2.96e-6 & 7.67e-6 & 65.7 & 5.12e-6 & 1.20e-5 & 41.9 & 0.64\\\cline{4-11}

   &      &     &  pwht & 5.11e-6 & 7.98e-6 & 87.9 & 7.07e-6 & 1.04e-5 & 60.1 & 0.68\\\cline{2-11}

   &   5\%   &   $(40\%,6.70)$  &  pdct & 1.07e-5 & 4.04e-5 & 410.0 & 1.00e-5 & 3.66e-5 & 320.1 & 0.78\\\cline{4-11}

   &      &     &  pfft & 1.34e-5 & 4.60e-5 & 351.8 & 1.10e-5 & 3.72e-5 & 279.5 & 0.79\\\cline{4-11}

   &      &     &  pwht & 1.03e-5 & 4.20e-5 & 411.2 & 1.04e-5 & 3.78e-5 & 322.6 & 0.78\\\cline{3-11}

   &      &   $(60\%,10.04)$  &  pdct & 8.08e-6 & 2.93e-5 & 191.1 & 5.29e-6 & 1.95e-5 & 150.0 & 0.78\\\cline{4-11}

   &      &     &  pfft & 7.71e-6 & 3.29e-5 & 187.9 & 5.05e-6 & 1.92e-5 & 136.7 & 0.73\\\cline{4-11}

   &      &     &  pwht & 8.58e-6 & 2.60e-5 & 204.3 & 5.26e-6 & 2.01e-5 & 150.7 & 0.74\\\cline{3-11}

   &      &   $(80\%,13.39)$  &  pdct & 6.02e-6 & 1.51e-5 & 106.5 & 3.67e-6 & 1.49e-5 & 77.4 & 0.73\\\cline{4-11}

   &      &     &  pfft & 5.69e-6 & 1.58e-5 & 100.0 & 3.51e-6 & 1.50e-5 & 70.2 & 0.70\\\cline{4-11}

   &      &     &  pwht & 5.34e-6 & 1.02e-5 & 107.6 & 3.96e-6 & 1.53e-5 & 77.1 & 0.72\\\cline{2-11}

   &   10\%   &   $(40\%,3.64)$  &  pdct & 2.72e-2 & 7.56e-2 & \mc{1}{|c||}{---} & 1.21e-5 & 4.04e-5 & 886.7 & ---\\\cline{4-11}

   &      &     &  pfft & 1.53e-5 & 5.22e-5 & 756.4 & 1.22e-5 & 4.12e-5 & 562.9 & 0.74\\\cline{4-11}

   &      &     &  pwht & 2.32e-2 & 6.44e-2 & \mc{1}{|c||}{---} & 1.22e-5 & 4.02e-5 & 863.8 & ---\\\cline{3-11}

   &      &   $(60\%,5.47)$  &  pdct & 8.47e-6 & 2.93e-5 & 345.7 & 7.20e-6 & 2.55e-5 & 243.7 & 0.70\\\cline{4-11}

   &      &     &  pfft & 9.05e-6 & 3.35e-5 & 310.9 & 7.43e-6 & 2.38e-5 & 240.9 & 0.77\\\cline{4-11}

   &      &     &  pwht & 8.74e-6 & 3.07e-5 & 345.0 & 7.90e-6 & 2.47e-5 & 243.5 & 0.71\\\cline{3-11}

   &      &   $(80\%,7.29)$  &  pdct & 7.99e-6 & 1.85e-5 & 170.7 & 6.71e-6 & 1.91e-5 & 123.0 & 0.72\\\cline{4-11}

   &      &     &  pfft & 8.08e-6 & 2.02e-5 & 160.3 & 6.17e-6 & 1.52e-5 & 115.6 & 0.72\\\cline{4-11}

   &      &     &  pwht & 8.00e-6 & 1.88e-5 & 170.5 & 6.69e-6 & 1.85e-5 & 123.2 & 0.72\\\cline{1-11}

%% file: results/18-Jul-2014/r_10.dat
$10$   &   1\%   &   $(40\%,13.59)$  &  pdct & 2.26e-5 & 5.31e-5 & 341.8 & 1.96e-5 & 5.85e-5 & 250.9 & 0.73\\\cline{4-11}

   &      &     &  pfft & 2.14e-5 & 4.89e-5 & 331.2 & 2.37e-5 & 6.37e-5 & 237.6 & 0.72\\\cline{4-11}

   &      &     &  pwht & 2.49e-5 & 6.29e-5 & 337.1 & 1.81e-5 & 5.52e-5 & 252.5 & 0.75\\\cline{3-11}

   &      &   $(60\%,20.38)$  &  pdct & 2.64e-5 & 4.49e-5 & 203.3 & 1.39e-5 & 3.33e-5 & 145.8 & 0.72\\\cline{4-11}

   &      &     &  pfft & 1.28e-5 & 3.28e-5 & 159.2 & 8.82e-6 & 2.90e-5 & 115.1 & 0.72\\\cline{4-11}

   &      &     &  pwht & 3.02e-5 & 5.61e-5 & 194.7 & 1.28e-5 & 3.28e-5 & 139.5 & 0.72\\\cline{3-11}

   &      &   $(80\%,27.18)$  &  pdct & 3.16e-6 & 1.05e-5 & 107.9 & 1.15e-5 & 2.56e-5 & 68.9 & 0.64\\\cline{4-11}

   &      &     &  pfft & 4.21e-6 & 1.35e-5 & 77.7 & 7.08e-6 & 1.93e-5 & 45.3 & 0.58\\\cline{4-11}

   &      &     &  pwht & 1.56e-5 & 2.61e-5 & 94.6 & 6.58e-6 & 1.43e-5 & 62.8 & 0.66\\\cline{2-11}

   &   5\%   &   $(40\%,5.76)$  &  pdct & 8.98e-6 & 3.71e-5 & 540.6 & 1.40e-5 & 5.08e-5 & 409.7 & 0.76\\\cline{4-11}

   &      &     &  pfft & 1.35e-5 & 5.30e-5 & 448.1 & 1.16e-5 & 4.38e-5 & 346.8 & 0.77\\\cline{4-11}

   &      &     &  pwht & 9.14e-6 & 3.75e-5 & 535.9 & 1.41e-5 & 4.92e-5 & 406.9 & 0.76\\\cline{3-11}

   &      &   $(60\%,8.64)$  &  pdct & 7.32e-6 & 3.32e-5 & 240.2 & 9.87e-6 & 3.38e-5 & 182.9 & 0.76\\\cline{4-11}

   &      &     &  pfft & 8.58e-6 & 3.50e-5 & 222.4 & 1.02e-5 & 3.43e-5 & 170.0 & 0.76\\\cline{4-11}

   &      &     &  pwht & 3.60e-5 & 5.48e-5 & 251.9 & 7.68e-6 & 2.67e-5 & 187.8 & 0.75\\\cline{3-11}

   &      &   $(80\%,11.52)$  &  pdct & 4.63e-6 & 1.31e-5 & 121.1 & 4.60e-6 & 1.86e-5 & 84.7 & 0.70\\\cline{4-11}

   &      &     &  pfft & 6.96e-6 & 2.84e-5 & 109.5 & 4.27e-6 & 1.82e-5 & 76.9 & 0.70\\\cline{4-11}

   &      &     &  pwht & 3.62e-6 & 1.15e-5 & 122.0 & 4.45e-6 & 1.87e-5 & 84.7 & 0.69\\\cline{2-11}

   &   10\%   &   $(40\%,3.35)$  &  pdct & 7.05e-2 & 2.31e-1 & \mc{1}{|c||}{---} & 2.74e-2 & 9.34e-2 & \mc{1}{|c||}{---} & ---\\\cline{4-11}

   &      &     &  pfft & 1.10e-5 & 4.21e-5 & 967.3 & 1.29e-5 & 4.45e-5 & 703.9 & 0.73\\\cline{4-11}

   &      &     &  pwht & 6.75e-2 & 2.19e-1 & \mc{1}{|c||}{---} & 2.43e-2 & 8.20e-2 & \mc{1}{|c||}{---} & ---\\\cline{3-11}

   &      &   $(60\%,5.02)$  &  pdct & 7.33e-6 & 3.27e-5 & 399.3 & 7.76e-6 & 2.77e-5 & 301.0 & 0.75\\\cline{4-11}

   &      &     &  pfft & 9.73e-6 & 3.59e-5 & 353.5 & 7.65e-6 & 2.73e-5 & 267.2 & 0.76\\\cline{4-11}

   &      &     &  pwht & 8.36e-6 & 3.24e-5 & 397.6 & 7.56e-6 & 2.72e-5 & 298.9 & 0.75\\\cline{3-11}

   &      &   $(80\%,6.70)$  &  pdct & 8.13e-6 & 2.03e-5 & 186.1 & 6.94e-6 & 2.25e-5 & 132.2 & 0.71\\\cline{4-11}

   &      &     &  pfft & 8.63e-6 & 3.00e-5 & 172.9 & 6.83e-6 & 2.02e-5 & 121.5 & 0.70\\\cline{4-11}

   &      &     &  pwht & 8.67e-6 & 2.43e-5 & 185.1 & 7.39e-6 & 1.89e-5 & 132.2 & 0.71\\\cline{1-11}

%% file: results/18-Jul-2014/r_15.dat
$15$   &   1\%   &   $(40\%,10.23)$  &  pdct & 1.76e-5 & 5.77e-5 & 406.3 & 1.47e-5 & 5.32e-5 & 309.7 & 0.76\\\cline{4-11}

   &      &     &  pfft & 2.35e-5 & 6.64e-5 & 374.7 & 1.55e-5 & 5.51e-5 & 279.0 & 0.74\\\cline{4-11}

   &      &     &  pwht & 1.59e-5 & 5.41e-5 & 401.0 & 1.51e-5 & 5.54e-5 & 305.3 & 0.76\\\cline{3-11}

   &      &   $(60\%,15.35)$  &  pdct & 2.96e-5 & 5.79e-5 & 225.6 & 1.48e-5 & 3.83e-5 & 159.5 & 0.71\\\cline{4-11}

   &      &     &  pfft & 9.59e-6 & 2.73e-5 & 188.7 & 9.38e-6 & 3.57e-5 & 134.3 & 0.71\\\cline{4-11}

   &      &     &  pwht & 2.76e-5 & 5.05e-5 & 238.8 & 1.40e-5 & 3.43e-5 & 168.7 & 0.71\\\cline{3-11}

   &      &   $(80\%,20.47)$  &  pdct & 1.58e-5 & 2.77e-5 & 123.7 & 8.52e-6 & 1.32e-5 & 80.4 & 0.65\\\cline{4-11}

   &      &     &  pfft & 9.77e-6 & 1.36e-5 & 94.4 & 8.31e-6 & 1.74e-5 & 53.9 & 0.57\\\cline{4-11}

   &      &     &  pwht & 1.49e-5 & 2.36e-5 & 117.8 & 1.23e-5 & 2.68e-5 & 75.7 & 0.64\\\cline{2-11}

   &   5\%   &   $(40\%,5.06)$  &  pdct & 2.10e-5 & 7.95e-5 & 674.6 & 1.34e-5 & 5.11e-5 & 506.3 & 0.75\\\cline{4-11}

   &      &     &  pfft & 1.30e-5 & 5.44e-5 & 546.1 & 1.68e-5 & 5.96e-5 & 412.9 & 0.76\\\cline{4-11}

   &      &     &  pwht & 2.07e-5 & 7.39e-5 & 674.5 & 1.32e-5 & 5.13e-5 & 505.5 & 0.75\\\cline{3-11}

   &      &   $(60\%,7.59)$  &  pdct & 7.78e-6 & 3.51e-5 & 280.9 & 1.04e-5 & 3.61e-5 & 208.2 & 0.74\\\cline{4-11}

   &      &     &  pfft & 9.38e-6 & 3.94e-5 & 254.7 & 1.12e-5 & 3.63e-5 & 190.9 & 0.75\\\cline{4-11}

   &      &     &  pwht & 6.72e-6 & 3.39e-5 & 283.1 & 1.03e-5 & 3.73e-5 & 210.9 & 0.75\\\cline{3-11}

   &      &   $(80\%,10.12)$  &  pdct & 5.39e-6 & 1.48e-5 & 135.8 & 5.15e-6 & 1.96e-5 & 93.1 & 0.69\\\cline{4-11}

   &      &     &  pfft & 7.64e-6 & 2.68e-5 & 120.7 & 5.61e-6 & 1.83e-5 & 82.9 & 0.69\\\cline{4-11}

   &      &     &  pwht & 6.62e-6 & 2.03e-5 & 134.7 & 4.98e-6 & 2.14e-5 & 92.9 & 0.69\\\cline{2-11}

   &   10\%   &   $(40\%,3.10)$  &  pdct & 1.02e-1 & 3.52e-1 & \mc{1}{|c||}{---} & 6.25e-2 & 2.27e-1 & \mc{1}{|c||}{---} & ---\\\cline{4-11}

   &      &     &  pfft & 2.46e-2 & 7.71e-2 & \mc{1}{|c||}{---} & 1.43e-5 & 5.13e-5 & 879.9 & ---\\\cline{4-11}

   &      &     &  pwht & 1.00e-1 & 3.50e-1 & \mc{1}{|c||}{---} & 6.03e-2 & 2.23e-1 & \mc{1}{|c||}{---} & ---\\\cline{3-11}

   &      &   $(60\%,4.65)$  &  pdct & 7.79e-6 & 3.35e-5 & 452.6 & 8.05e-6 & 2.76e-5 & 334.7 & 0.74\\\cline{4-11}

   &      &     &  pfft & 9.62e-6 & 3.85e-5 & 395.7 & 8.35e-6 & 3.03e-5 & 295.3 & 0.75\\\cline{4-11}

   &      &     &  pwht & 7.97e-6 & 3.40e-5 & 455.6 & 7.76e-6 & 2.85e-5 & 336.3 & 0.74\\\cline{3-11}

   &      &   $(80\%,6.20)$  &  pdct & 8.61e-6 & 2.25e-5 & 203.2 & 7.53e-6 & 1.93e-5 & 143.3 & 0.70\\\cline{4-11}

   &      &     &  pfft & 9.44e-6 & 3.26e-5 & 186.0 & 7.78e-6 & 2.65e-5 & 128.1 & 0.69\\\cline{4-11}

   &      &     &  pwht & 9.02e-6 & 2.61e-5 & 202.1 & 8.06e-6 & 2.60e-5 & 142.6 & 0.71\\\cline{1-11}

%% file: results/18-Jul-2014/r_20.dat
$20$   &   1\%   &   $(40\%,8.22)$  &  pdct & 2.14e-5 & 5.94e-5 & 503.4 & 1.91e-5 & 6.27e-5 & 374.3 & 0.74\\\cline{4-11}

   &      &     &  pfft & 1.80e-5 & 6.55e-5 & 433.6 & 1.85e-5 & 6.98e-5 & 327.2 & 0.75\\\cline{4-11}

   &      &     &  pwht & 1.67e-5 & 5.88e-5 & 496.4 & 1.70e-5 & 6.27e-5 & 371.8 & 0.75\\\cline{3-11}

   &      &   $(60\%,12.33)$  &  pdct & 2.62e-5 & 6.27e-5 & 258.1 & 1.35e-5 & 4.07e-5 & 181.2 & 0.70\\\cline{4-11}

   &      &     &  pfft & 2.44e-5 & 5.63e-5 & 222.2 & 1.00e-5 & 3.75e-5 & 158.6 & 0.71\\\cline{4-11}

   &      &     &  pwht & 2.61e-5 & 5.51e-5 & 260.6 & 2.04e-5 & 5.05e-5 & 179.8 & 0.69\\\cline{3-11}

   &      &   $(80\%,16.43)$  &  pdct & 1.36e-5 & 3.26e-5 & 126.6 & 8.74e-6 & 2.52e-5 & 80.5 & 0.64\\\cline{4-11}

   &      &     &  pfft & 1.04e-5 & 2.38e-5 & 98.2 & 5.29e-6 & 2.18e-5 & 55.2 & 0.56\\\cline{4-11}

   &      &     &  pwht & 1.41e-5 & 3.16e-5 & 116.1 & 5.37e-6 & 2.23e-5 & 74.7 & 0.64\\\cline{2-11}

   &   5\%   &   $(40\%,4.51)$  &  pdct & 1.90e-5 & 7.53e-5 & 835.2 & 1.70e-5 & 6.54e-5 & 614.5 & 0.74\\\cline{4-11}

   &      &     &  pfft & 1.25e-5 & 5.42e-5 & 654.1 & 1.74e-5 & 6.67e-5 & 487.6 & 0.75\\\cline{4-11}

   &      &     &  pwht & 1.86e-5 & 7.64e-5 & 830.1 & 1.73e-5 & 6.55e-5 & 611.9 & 0.74\\\cline{3-11}

   &      &   $(60\%,6.77)$  &  pdct & 7.43e-6 & 3.62e-5 & 326.8 & 1.06e-5 & 4.26e-5 & 238.7 & 0.73\\\cline{4-11}

   &      &     &  pfft & 9.77e-6 & 4.27e-5 & 291.6 & 1.15e-5 & 3.87e-5 & 214.0 & 0.73\\\cline{4-11}

   &      &     &  pwht & 7.75e-6 & 3.68e-5 & 326.4 & 1.09e-5 & 4.34e-5 & 238.3 & 0.73\\\cline{3-11}

   &      &   $(80\%,9.02)$  &  pdct & 7.37e-6 & 2.65e-5 & 147.8 & 5.48e-6 & 2.25e-5 & 101.6 & 0.69\\\cline{4-11}

   &      &     &  pfft & 7.63e-6 & 3.67e-5 & 132.7 & 5.28e-6 & 2.52e-5 & 88.9 & 0.67\\\cline{4-11}

   &      &     &  pwht & 7.76e-6 & 3.13e-5 & 147.9 & 5.25e-6 & 2.37e-5 & 101.8 & 0.69\\\cline{2-11}

   &   10\%   &   $(40\%,2.88)$  &  pdct & 1.32e-1 & 4.68e-1 & \mc{1}{|c||}{---} & 9.38e-2 & 3.55e-1 & \mc{1}{|c||}{---} & ---\\\cline{4-11}

   &      &     &  pfft & 5.97e-2 & 1.96e-1 & \mc{1}{|c||}{---} & 1.20e-2 & 4.18e-2 & \mc{1}{|c||}{---} & ---\\\cline{4-11}

   &      &     &  pwht & 1.32e-1 & 4.65e-1 & \mc{1}{|c||}{---} & 9.32e-2 & 3.51e-1 & \mc{1}{|c||}{---} & ---\\\cline{3-11}

   &      &   $(60\%,4.33)$  &  pdct & 6.99e-6 & 3.16e-5 & 517.3 & 1.20e-5 & 4.24e-5 & 375.2 & 0.73\\\cline{4-11}

   &      &     &  pfft & 9.84e-6 & 4.06e-5 & 441.5 & 8.81e-6 & 3.23e-5 & 325.2 & 0.74\\\cline{4-11}

   &      &     &  pwht & 7.16e-6 & 3.22e-5 & 512.4 & 1.25e-5 & 4.18e-5 & 372.7 & 0.73\\\cline{3-11}

   &      &   $(80\%,5.77)$  &  pdct & 9.26e-6 & 2.73e-5 & 219.0 & 8.39e-6 & 2.35e-5 & 153.3 & 0.70\\\cline{4-11}

   &      &     &  pfft & 9.77e-6 & 3.32e-5 & 200.4 & 8.07e-6 & 2.81e-5 & 137.5 & 0.69\\\cline{4-11}

   &      &     &  pwht & 9.13e-6 & 3.38e-5 & 220.6 & 8.17e-6 & 2.23e-5 & 154.9 & 0.70\\\cline{1-11}